\documentclass{article}

\usepackage{amsmath,amssymb,amsfonts,theorem}
\usepackage{graphicx,color}
\usepackage{subfigure}
\graphicspath{{./Figures/}}
\usepackage{graphicx}
\usepackage{color}
\usepackage{amsfonts}
\usepackage{amssymb}
\usepackage{gastex}

{ \theorembodyfont{\normalfont} 

\newtheorem{remark}{Remark}
}

\newtheorem{definition}{Definition}

\newtheorem{lemma}{Lemma}

\newtheorem{proposition}{Proposition}

\newcommand{\reals}{\mathbb{R}} 

\newcommand{\real}{\ensuremath{\mathbb{R}}}

\newcommand{\R}{\reals}

\newcommand{\dst}{\displaystyle}

\def\be{\begin{equation}}
\def\ee{\end{equation}}
\def\ba{\begin{array}}
\def\ea{\end{array}}
\def\eqa{\begin{eqnarray}}
\def\eqe{\end{eqnarray}}

\author{ \parbox{4 in}{\centering Claudio De Persis
           and Fr\'{e}d\'{e}ric Mazenc
          \thanks{C. De Persis is with the Faculty of Engineering Technology (CTW), University of Twente, the Netherlands ({\tt\small
         c.depersis@ctw.utwente.nl}) and with the Department of Computer and System Sciences, Sapienza University of Rome, Italy.
         F.\ Mazenc is with Projet INRIA DISCO, CNRS-Sup\'elec, 3 rue Joliot Curie,
91192 Gif-sur-Yvette, France ({\tt\small
Frederic.Mazenc@lss.supelec.fr}).}
          } }

\title{\LARGE \bf
Stability of quantized time-delay nonlinear systems: A
Lyapunov-Krasowskii-functional approach}

\begin{document}

\maketitle


\begin{abstract}
Lyapunov-Krasowskii functionals are used to design quantized
control laws for nonlinear continuous-time systems in the presence
of constant delays in the input. The quantized control law is
implemented via hysteresis to avoid chattering. Under appropriate
conditions, our analysis applies to  stabilizable nonlinear
systems for any value of the quantization density. The resulting
quantized feedback is parametrized with respect to the
quantization density. Moreover, the maximal allowable delay
tolerated by the system is characterized as a function of the
quantization density.
\end{abstract}

\medskip


\section{Introduction}
Quantized control systems (\cite{elia.mitter.tac01},
\cite{liberzon.aut03}), are systems in which the control law is a
piece-wise constant function of time taking values in a finite
set. The design of quantized control systems is based on a
partition of the state space. One value of the control law is
associated to each set of the partition, and whenever the state
crosses the boundary between two sets of the partition, the
control law takes the new
value associated to the set which the state has just entered.\\
When dealing with the problem of stabilizing the origin of the
state space for linear discrete-time systems, the paper
\cite{elia.mitter.tac01} has shown the effectiveness of
logarithmic quantization in which the partition of the state space
is coarser away from the origin and denser in its vicinity. It has
also introduced the notion of quantization density, that is the
number of regions of the partition per unit of space. Intuitively,
the larger is the quantization density, the easier is the
quantized control problem, since as the quantization density gets
larger, the quantized control law approaches a control law without
quantization. The paper \cite{liberzon.aut03} deals with a similar
problem but for nonlinear continuous-time systems which can be
made input-to-state stable with respect to the quantization error.
Recently, the paper \cite{ceragioli.depersis.scl07} has
investigated quantized control systems in the framework of
discontinuous control systems, discussing appropriate notions of
solutions, namely Krasowskii and Carath\'eodory solutions. In this
framework, the effect of quantization is viewed as an additional
disturbance whose effect is attenuated by a Lyapunov redesign of
the control law. Namely, given any nonlinear continuous-time
process which is stabilizable by a continuous feedback, and given
any value of the quantization density, it is always possible to
find a new feedback depending on the quantization density, in such
a way that the process in closed-loop with the quantized control
law is practically stable with a basin of attraction which can be
made arbitrarily large. Other notions of robustness (namely,
robustness in the sense of the ${\cal L}_2$-gain) in connection
with quantized control problems have been examined in
\cite{hayakawa.et.al.acc06} and \cite{ceragioli.depersis.scl07}.
Moreover, in the former, an
adaptive quantized control scheme has been investigated.\\
Since quantized controls take values in a finite set, they lend
themselves to be implemented over a finite data-rate communication
channel. Data transmitted over a channel are usually delivered at
the other end of channel after a delay. The problem of quantized
control systems in the presence of delays then arises very
naturally. Such a problem has been examined for the first time in
\cite{liberzon.tac06}, where the connection between
Razumikhin-type theorems and the ISS small-gain theorem
established in \cite{teel.tac98} was exploited. In recent years,
besides \cite{teel.tac98}, other contributions in the area of
nonlinear time-delay systems have appeared (see, for instance,
\cite{michiels.et.al.sicon01}, \cite{sil1},
\cite{mazenc.et.al.tac04}, \cite{pepej}, \cite{jank}, \cite{kara},
\cite{mazenc.bliman.tac06}, \cite{fridman.et.al.aut08} and
references therein). In particular, the paper
\cite{mazenc.bliman.tac06} has proposed a
Lyapunov-Krasowskii-functional approach to study the
stabilizability of nonlinear systems in the presence of a delay in
the input.

The aim of this paper is to pursue the approach of
\cite{mazenc.bliman.tac06} in the analysis and design of quantized
time-delay control systems. Besides the use of Lyapunov-Krasowskii
functionals, there are other important features of the approach
which make our paper different from other contributions. We
implement the quantized control with the hysteretic mechanism
suggested in \cite{hayakawa.et.al.acc06} to avoid chattering. It
is known from \cite{depersis.scl09} that, in the case no delay is
present, the analysis of such hysteretic solutions can be reduced
to the analysis of Krasowskii and Carath\'eodory solutions
considered in \cite{ceragioli.depersis.scl07}. In the case of
quantized time-delay systems, the adoption of the hysteretic
solution is desirable. First, because it allows us to avoid
technical issues related to more general notions of solutions of
time-delay quantized (that is, discontinuous) systems. Second, the
existence of more general solutions such as Carath\'eodory
solutions is guaranteed only under additional conditions (see
e.g.\ \cite{ceragioli.depersis.scl07}). Another feature which is
worth mentioning is that, as in \cite{ceragioli.depersis.scl07},
our analysis applies to stabilizable nonlinear systems for any
value of the quantization density, provided that suitable
conditions are satisfied. Then, the quantized feedback which
stabilizes the closed-loop system despite the delay turns out to
be
parametrized with respect to the quantization density.\\
Our approach leads to a set of conditions to design quantized
control systems which are robust with respect to delays. Since we
employ the results of \cite{mazenc.bliman.tac06} based on
Lyapunov-Krasowskii functionals, our conditions represent an
alternative to the conditions derived using Razumikhin-like
theorems in \cite{teel.tac98}, \cite{liberzon.tac06}. Other
conditions could be derived using recent results on input-to-state
stability of time-delay systems via Lyapunov-Krasowskii
functionals (\cite{pepej}, and \cite{fridman.et.al.aut08} where a
few comments in this regard have been presented). However, this
investigation is beyond the scope of the paper.

In the next section, we present a few preliminaries, such as the
definition of the quantizer and the notion of solution we adopt.
The main result along with the standing assumptions and a couple
of examples are examined in Section 3. Proof of the main result is
given in Section 4. Conclusions are drawn in Section 5.

\medskip

\noindent {\bf Notation, definitions}
\begin{itemize}
\item $\R_{\ge 0}$ (respectively, $\R_{> 0}$) denotes the set of
non-negative (positive) real numbers.
\item Let $r_1,r_2$ be two real numbers such that $r_1 < r_2$.
Let $C^1([r_1,r_2], \R^m)$ (respectively, $\overline
C^1([r_1,r_2],$ $\R^m)$) denote the set of continuously
differentiable (respectively, piece-wise continuously
differentiable) functions $\phi(\cdot):[r_1,r_2]\to \R^m$.
\item {\em Norms}. $|\cdot|$ stands for the Euclidean norm,
$||\phi||_c=\sup_{t\in [r_1,r_2]} |\phi(t)|$ stands for the norm
of a function $\phi\in C^1([r_1,r_2], \R^m)$.
\item ${\rm sgn}(r)$, $r\in \mathbb{R}$, denotes the sign function, i.e.
the function such that ${\rm sgn}(r) = 1$ if $r > 0$, ${\rm
sgn}(r)= - 1$ if $r < 0$, and ${\rm sgn}(r) = 0$ if $r = 0$.
\item To simplify the notation we will frequently use the
notation of the Lie derivative. More precisely, if $f :
\mathbb{R}^n \rightarrow \mathbb{R}^n$ is a vector field and $h :
\mathbb{R}^n \rightarrow \mathbb{R}$ is a scalar function, we may
use the notation $L_f h(x)$ for $\frac{\partial h}{\partial x}(x)
f(x)$.
\item A continuous function $k:[0,\infty)\to [0,\infty)$ is of {\em class ${\cal K}$}
provided it is  zero at zero and strictly increasing. A {\em class
${\cal K}_{\infty}$} function is a {\em class ${\cal K}$} function
which in addition is  unbounded.
\item We shall often omit arguments of functions to simplify notation.
\item For a real-valued function $z(t)$, we denote by $z(t^+)$
the right limit\linebreak $\lim_{m>t,m\to t} z(m)$.
\end{itemize}

\section{Problem formulation}
\label{secp}

We are interested in investigating the stability property of
systems when the feedback control law undergoes quantization and
delays. This problem arises in (idealized) scenarios in which a
finite bandwidth channel lies in the feedback loop and introduces
a delay. In the sub-sections below, we recall what is meant by
quantization and what is a quantizer, we introduce the quantized
time-delay system and the notion of solution we adopt, and finally
the formulation of the problem.
\subsection{Quantizers}
To the purpose of describing our system in more formal terms, we
introduce  the following multi-valued map, which will be referred
to henceforth as the {\em quantizer}. Let $u_0>0$ and $ 0<\rho<1$
be real numbers, let $u_i=\rho^i u_0$ and $U=\{ 0, \pm u_i, \pm
u_i(1+\delta)^{-1},$ $i=0,1,\ldots,j\}$, with $j\ge 1$ an integer.
Let $\delta=(1-\rho)(1+\rho)^{-1}$ and
\be\label{ppsi} \Psi
(u) = \left\{\ba{lrclllll} u_i{\rm sgn}(u),
&\dst\frac{1}{1+\delta}u_i& < & |u| &\le &
\dst\frac{1}{1-\delta}u_i\;, & 0\le i\le j \\[2mm]
\dst\frac{u_i}{1+\delta}{\rm sgn}(u),
&\dst\frac{1}{(1+\delta)^2}u_i& < & |u| &\le &
\dst\frac{1}{(1+\delta)(1-\delta)}u_i\;, & 0\le i\le j \\[2mm]
0, & 0 &\le & |u| &\le & \dst\frac{1}{1+\delta}u_j\;. & \\[2mm]
\ea\right. \ee A picture of the map is given in Fig.\
\ref{fig.multivalued}. Observe for later use that \be\label{aeio}
\rho = \frac{1 - \delta}{1 + \delta} \ee and \be\label{yeio} u_i =
\left(\frac{1 - \delta}{1 + \delta}\right)^i u_0 \; , \; \forall i
\in \{0,1,\ldots,j \} \; . \ee A few remarks are in order:
\begin{itemize}
\item The {\em range} of the quantizer, i.e.\ its interval of definition, is
$[-\frac{u_0}{1-\delta}, \frac{u_0}{1-\delta}]$. We do not define
$\Psi(u)$ for $|u|> \frac{u_0}{1-\delta}$, since we will design
the parameter $u_0$ in such a way that the control $|u(t-\tau)|$,
which is the actual argument of the map $\Psi$, never exceeds this
upper bound.
\item The logarithmic quantizer with a {\em finite} number of quantization
levels, which is a truncated version of the quantizer with an
infinite number of quantization levels, was introduced  in
\cite{elia.mitter.tac01}, Section V, and it is as follows:
\be\label{ppsi.elia.mitter} \Psi(u) = \left\{\ba{lrclrcll} u_i{\rm
sgn}(u), &\dst\frac{1}{1+\delta}u_i& < & |u| &\le &
\dst\frac{1}{1-\delta}u_i, & 0\le i\le j \\[2mm]
0, & 0 &\le & |u| &\le & \dst\frac{1}{1+\delta}u_j\;. & \\[2mm]
\ea\right. \ee Compared with (\ref{ppsi.elia.mitter}), the
quantizer (\ref{ppsi}) considered in this paper has additional
quantization levels. To have a pictorial representation of the
quantizer (\ref{ppsi.elia.mitter}), one can refer to Fig.\
\ref{fig.multivalued} and remove the quantization levels labeled
as $\frac{u_0}{1+\delta}$ and $\frac{u_1}{1+\delta}$. The new
quantization levels in (\ref{ppsi}) are added to avoid chattering.
This will be explained in detail as soon as  the system we are
interested in and the notion of solution we adopt are introduced
(see Remark \ref{rm1} below).
\item The parameter $\rho$ can be viewed as a measure of the quantization {\em density},
since the smaller is $\rho$,  the coarser is the quantizer
(\cite{elia.mitter.tac01}). In fact, by (\ref{aeio}), as $\rho$
approaches $0$, $\delta$ approaches $1$, that is the width of the
sector bound in Fig.\  \ref{fig.multivalued} gets larger and,
given an interval of fixed length on the $u$-axis in Fig.\
\ref{fig.multivalued}, $\Psi(u)$ will have fewer quantization
levels as $u$ ranges over that interval.
\item In the quantizer (\ref{ppsi}), the parameters $\delta, u_0,j$
appear. Throughout the paper, we shall assume that $\delta$ can
take any value in the interval $(0,1)$ (i.e.\ the quantization
density can be equal to any value). On the other hand the positive
real number $u_0$ (which defines the range of the quantizer) and
the integer $j$ (which gives the number of quantization levels)
are to be designed. Although it would be more correct to denote
explicitly the dependence of $\Psi$ on $u_0, j$, i.e.\ to have
$\Psi_{j,u_0}(u)$, this is not pursued in the paper to avoid
cumbersome notations.
\end{itemize}

\subsection{Quantized time-delay systems}

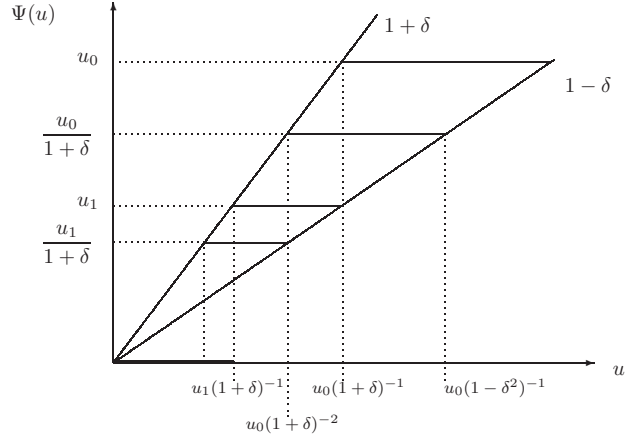
\begin{figure}
\setlength{\unitlength}{1mm}
\begin{center}
\scalebox{0.8}{
\begin{picture}(70,80)(0,0)
\put(10,10){\vector(1,0){80}} \put(10,10){\vector(0,1){60}}
\put(-7,67){$\Psi (u)$} \put(93,8){$u$}
%

\multiput(10,10)(.049,.03373819163){1495}{\line(1,0){.04723346829}}
\put(85,55){$1-\delta$}
\multiput(10,10)(.03372891943,.0445){1300}{\line(0,1){.04372267333}}
\put(55,65){$1+\delta$}

\put(48,60){\line(1,0){34.5}}
\multiput(48,60)(-1,0){38}{{\rule{.4pt}{.4pt}}} \put(4,60){$u_0$}

\put(39,48){\line(1,0){26}}
\multiput(39,48)(-1,0){30}{{\rule{.4pt}{.4pt}}}
\put(-2,47){$\dst\frac{u_0}{1+\delta}$}

\put(25,30){\line(1,0){14}}
\multiput(10,30)(1,0){16}{{\rule{.4pt}{.4pt}}}
\put(-2,29){$\dst\frac{u_{1}}{1+\delta}$}

\multiput(25,29)(0,-1){20}{{\rule{.4pt}{.4pt}}}

\multiput(39,48)(0,-1){48}{{\rule{.4pt}{.4pt}}}
\put(32,-1){{\footnotesize $u_0(1+\delta)^{-2}$}}

\put(47.9,36.1){\line(-1,0){18.25}}
\multiput(10,36)(1,0){19}{{\rule{.4pt}{.4pt}}} \put(4,36){$u_{1}$}

\put(10,10.1){\line(1,0){20}} \put(10,10.2){\line(1,0){20}}
\put(10,10.3){\line(1,0){20}}

\multiput(48,7)(0,1){53}{{\rule{.4pt}{.4pt}}}
\put(43,5){{\footnotesize $u_0(1+\delta)^{-1}$}}

\multiput(65,7)(0,1){42}{{\rule{.4pt}{.4pt}}}
\put(65,5){{\footnotesize $u_0(1-\delta^2)^{-1}$}}

\multiput(30,36)(0,-1){30}{{\rule{.4pt}{.4pt}}}
\put(23,5){{\footnotesize $u_{1}(1+\delta)^{-1}$}}

\end{picture}
} \caption{\label{fig.multivalued}The multi-valued map $\Psi (u)$
for $u>0$, and with $j=1$.}
\end{center}
\end{figure}

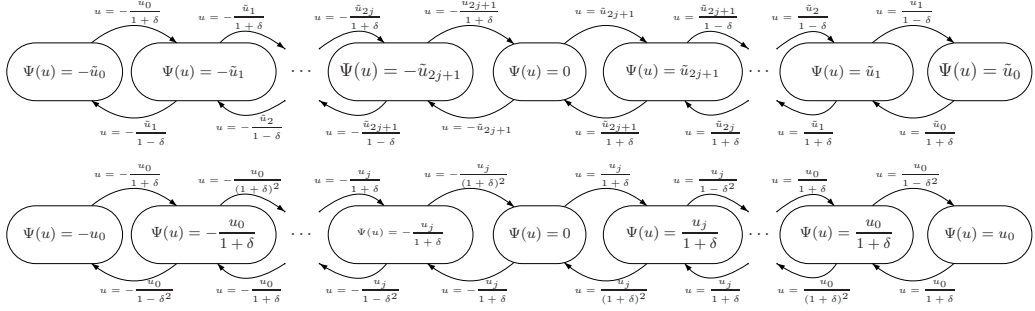
\begin{figure}
\vspace{1cm}
\setlength{\unitlength}{1.1mm}
\hspace{1.5cm}
\begin{tabular}{c}
\scalebox{0.7}{
\begin{picture}(88,20)(14,-62)


\node[NLangle=0.0,Nw=20.0,Nh=10.0,Nmr=10.0](n0)(2,-48.0){\footnotesize{$\Psi(u)=-\tilde
u_0$}}

\node[NLangle=0.0,Nw=25.0,Nh=10.0,Nmr=10.0](n3)(26.0,-48.0)
{\footnotesize{$\Psi(u)=-\tilde u_1$}}

\drawedge[curvedepth=8.0](n0,n3){}

\put(7,-38){\tiny{$u=-\dst\frac{\tilde u_0}{1+\delta}$}}

\drawedge[curvedepth=8.0](n3,n0){}

\put(8,-59){\tiny{$u=-\dst\frac{\tilde u_1}{1-\delta}$}}

\node[NLangle=0.0,Nw=0.1,Nh=0.1,Nmr=10.0](n43)(40.0,-44.0){}

\drawedge[curvedepth=5.0](n3,n43){}

\put(25,-39){\tiny{$u=-\dst\frac{\tilde u_1}{1+\delta}$}}

\node[NLangle=0.0,Nw=0.1,Nh=0.1,Nmr=10.0](n34)(40.0,-53.0){}

\drawedge[curvedepth=5.0](n34,n3){}

\put(28,-58){\tiny{$u=-\dst\frac{\tilde u_2}{1-\delta}$}}

\node[NLangle=0.0,Nw=0.1,Nh=0.1,Nmr=10.0](n10)(46.0,-44.0){}

\node[NLangle=0.0,Nw=0.1,Nh=0.1,Nmr=10.0](n11)(46.0,-53.0){}

\put(41,-48){$\ldots$}

\node[NLangle=0.0,Nw=25.0,Nh=10.0,Nmr=10.0](n4)(60.0,-48.0)
{$\Psi(u)=-\tilde u_{2j+1}$}

\drawedge[curvedepth=5.0](n10,n4){}
\put(45,-39){\tiny{$u=-\dst\frac{\tilde u_{2j}}{1+\delta}$}}

\drawedge[curvedepth=5.0](n4,n11){}
\put(47,-59){\tiny{$u=-\dst\frac{\tilde u_{2j+1}}{1-\delta}$}}

\node[NLangle=0.0,Nw=16.0,Nh=10.0,Nmr=10.0](n5)(84.0,-48.0)
{\footnotesize{$\Psi(u)=0$}}

\drawedge[curvedepth=8.0](n4,n5){}

\put(64,-38){\tiny{$u=-\dst\frac{\tilde u_{2j+1}}{1+\delta}$}}

\drawedge[curvedepth=8.0](n5,n4){}

\put(67,-58){\tiny{$u=-\tilde u_{2j+1}$}}

\node[NLangle=0.0,Nw=24.0,Nh=10.0,Nmr=10.0](n6)(107.0,-48.0)
{\footnotesize{$\Psi(u)=\tilde u_{2j+1}$}}

\node[NLangle=0.0,Nw=0.1,Nh=0.1,Nmr=10.0](n65)(120.0,-44.0){}
\drawedge[curvedepth=5.0](n6,n65){}
\put(107,-39){\tiny{$u=\dst\frac{\tilde u_{2j+1}}{1-\delta}$}}

\node[NLangle=0.0,Nw=0.1,Nh=0.1,Nmr=10.0](n56)(120.0,-53.0){}
\drawedge[curvedepth=5.0](n56,n6){}
\put(109,-59){\tiny{$u=\dst\frac{\tilde u_{2j}}{1+\delta}$}}

\drawedge[curvedepth=8.0](n5,n6){}

\put(90,-38){\tiny{$u=\tilde u_{2j+1}$}}

\drawedge[curvedepth=8.0](n6,n5){}

\put(90,-59){\tiny{$u=\dst\frac{\tilde u_{2j+1}}{1+\delta}$}}

\put(120,-48){$\ldots$}

\node[NLangle=0.0,Nw=23.0,Nh=10.0,Nmr=10.0](n7)(137.0,-48.0)
{\footnotesize{$\Psi(u)=\tilde u_1$}}

\node[NLangle=0.0,Nw=0.1,Nh=0.1,Nmr=10.0](n75)(125.0,-44.0){}
\drawedge[curvedepth=5.0](n75,n7){}
\put(124,-39){\tiny{$u=\dst\frac{\tilde u_2}{1-\delta}$}}

\node[NLangle=0.0,Nw=0.1,Nh=0.1,Nmr=10.0](n57)(125.0,-53.0){}
\drawedge[curvedepth=5.0](n7,n57){}
\put(125,-59){\tiny{$u=\dst\frac{\tilde u_1}{1+\delta}$}}

\node[NLangle=0.0,Nw=18.0,Nh=10.0,Nmr=10.0](n8)(160.0,-48.0)
{$\Psi(u)=\tilde u_0$}

\drawedge[curvedepth=8.0](n7,n8){}

\put(142,-38){\tiny{$u=\dst\frac{\tilde u_1}{1-\delta}$}}

\drawedge[curvedepth=8.0](n8,n7){}

\put(146,-59){\tiny{$u=\dst\frac{\tilde u_0}{1+\delta}$}}

\end{picture}
}
\\[5mm]
\scalebox{0.7}{
\begin{picture}(88,20)(14,-62)


\node[NLangle=0.0,Nw=20.0,Nh=10.0,Nmr=10.0](n0)(2,-48.0){\footnotesize{$\Psi(u)=-u_0$}}

\node[NLangle=0.0,Nw=25.0,Nh=10.0,Nmr=10.0](n3)(26.0,-48.0){\footnotesize{$\Psi(u)=-\dst\frac{u_0}{1+\delta}$}}

\drawedge[curvedepth=8.0](n0,n3){}

\put(7,-38){\tiny{$u=-\dst\frac{u_0}{1+\delta}$}}

\drawedge[curvedepth=8.0](n3,n0){}

\put(8,-58){\tiny{$u=-\dst\frac{u_0}{1-\delta^2}$}}

\node[NLangle=0.0,Nw=0.1,Nh=0.1,Nmr=10.0](n43)(40.0,-44.0){}

\drawedge[curvedepth=5.0](n3,n43){}

\put(25,-39){\tiny{$u=-\dst\frac{u_0}{(1+\delta)^2}$}}

\node[NLangle=0.0,Nw=0.1,Nh=0.1,Nmr=10.0](n34)(40.0,-53.0){}

\drawedge[curvedepth=5.0](n34,n3){}

\put(28,-58){\tiny{$u=-\dst\frac{u_0}{1+\delta}$}}

\node[NLangle=0.0,Nw=0.1,Nh=0.1,Nmr=10.0](n10)(46.0,-44.0){}

\node[NLangle=0.0,Nw=0.1,Nh=0.1,Nmr=10.0](n11)(46.0,-53.0){}

\put(41,-48){$\ldots$}

\node[NLangle=0.0,Nw=25.0,Nh=10.0,Nmr=10.0](n4)(60.0,-48.0)
{\tiny{$\Psi(u)=-\dst\frac{u_j}{1+\delta}$}}

\drawedge[curvedepth=5.0](n10,n4){}
\put(45,-39){\tiny{$u=-\dst\frac{u_{j}}{1+\delta}$}}

\drawedge[curvedepth=5.0](n4,n11){}
\put(47,-58){\tiny{$u=-\dst\frac{u_{j}}{1-\delta^2}$}}

\node[NLangle=0.0,Nw=16.0,Nh=10.0,Nmr=10.0](n5)(84.0,-48.0)
{\footnotesize{$\Psi(u)=0$}}

\drawedge[curvedepth=8.0](n4,n5){}

\put(64,-38){\tiny{$u=-\dst\frac{u_j}{(1+\delta)^2}$}}

\drawedge[curvedepth=8.0](n5,n4){}

\put(67,-58){\tiny{$u=-\dst\frac{u_j}{1+\delta}$}}

\node[NLangle=0.0,Nw=24.0,Nh=10.0,Nmr=10.0](n6)(107.0,-48.0)
{\footnotesize{$\Psi(u)=\dst\frac{u_j}{1+\delta}$}}

\node[NLangle=0.0,Nw=0.1,Nh=0.1,Nmr=10.0](n65)(120.0,-44.0){}
\drawedge[curvedepth=5.0](n6,n65){}
\put(107,-39){\tiny{$u=\dst\frac{u_j}{1-\delta^2}$}}

\node[NLangle=0.0,Nw=0.1,Nh=0.1,Nmr=10.0](n56)(120.0,-53.0){}
\drawedge[curvedepth=5.0](n56,n6){}
\put(109,-58){\tiny{$u=\dst\frac{u_j}{1+\delta}$}}

\drawedge[curvedepth=8.0](n5,n6){}

\put(90,-38){\tiny{$u=\dst\frac{u_j}{1+\delta}$}}

\drawedge[curvedepth=8.0](n6,n5){}

\put(90,-58){\tiny{$u=\dst\frac{u_j}{(1+\delta)^2}$}}

\put(120,-48){$\ldots$}

\node[NLangle=0.0,Nw=23.0,Nh=10.0,Nmr=10.0](n7)(137.0,-48.0)
{\footnotesize{$\Psi(u)=\dst\frac{u_0}{1+\delta}$}}

\node[NLangle=0.0,Nw=0.1,Nh=0.1,Nmr=10.0](n75)(125.0,-44.0){}
\drawedge[curvedepth=5.0](n75,n7){}
\put(124,-39){\tiny{$u=\dst\frac{u_0}{1+\delta}$}}

\node[NLangle=0.0,Nw=0.1,Nh=0.1,Nmr=10.0](n57)(125.0,-53.0){}
\drawedge[curvedepth=5.0](n7,n57){}
\put(125,-58){\tiny{$u=\dst\frac{u_0}{(1+\delta)^2}$}}

\node[NLangle=0.0,Nw=18.0,Nh=10.0,Nmr=10.0](n8)(160.0,-48.0)
{\footnotesize{$\Psi(u)=u_0$}}

\drawedge[curvedepth=8.0](n7,n8){}

\put(142,-38){\tiny{$u=\dst\frac{u_0}{1-\delta^2}$}}

\drawedge[curvedepth=8.0](n8,n7){}

\put(146,-58){\tiny{$u=\dst\frac{u_0}{1+\delta}$}}

\end{picture}
}
\end{tabular}
\caption{\label{fig.multivalued2} The graph at the top illustrates
the law (\ref{Psi.law}) which describes the evolution of
$\Psi(u(t))$ as $u(t)=\bar z(t)$ varies. Each edge connects two
nodes, and is labeled with the condition which triggers the
transition from the starting node to the destination node. The
graph at the bottom illustrates the same law but with nodes and
edges now labeled making use of  the original values $u$ rather
than $\tilde u$. }
\end{figure}

We are interested in investigating the stability
 of the quantized time-delay system
\be\label{system.ol} \dot x(t) = f(x(t)) + g(x(t))
\Psi(u(t-\tau))\;, \ee with $x(t)\in \R^n$, $n\ge 1$, $f(x),g(x)$
locally Lipschitz functions, and $\tau$ a positive real number,
when $u(t) = z(x(t))$, with $z(\cdot)$ a continuously
differentiable real-valued function to be designed. Since
$\Psi(u(t-\tau))$ is a multi-valued function, we must specify the
rule by which $\Psi(u(t-\tau))$ takes value in $U$ depending on
its argument $u(t-\tau)$.


Consider the initial condition $\varphi\in C^1([-2\tau,0],\R^n)$
and let $T<\tau$ be a suitable positive number. For $t\in [0,T)$
we focus our attention on $\Psi(\bar z(t))$, where to ease the
notation we have set $\bar z(t) := z(\varphi(t-\tau))$. At time
$t=0$, depending on $|\bar z(0)|$, the value taken by the
quantizer is specified as follows: \be\label{ppsi0} \Psi(\bar
z(0)) = \left\{\ba{lrclrclr} u_i{\rm sgn}(\bar z(0)), &
\dst\frac{1}{1+\delta}u_i& < & |\bar z(0)| & \le &
\dst\frac{1}{1-\delta}u_i\;, & 0\le i\le j \\[2mm]
0, & 0&\le & |\bar z(0)| &\le & \dst\frac{1}{1+\delta}u_j\;. & \\[2mm]
\ea\right. \ee For all $t\in [0,T)$, we describe the law according
to which $\Psi(\bar z(t))$ evolves as the argument $\bar z(t)$
varies. Before that, in order to have a concise description, we
rename the quantization levels as follows:
\[
\tilde u_k:=\left\{ \ba{lrr}
u_{k/2} & k\;{\rm even}\\[3mm]
\dst\frac{u_{(k-1)/2}}{1+\delta} & k\;{\rm odd}, & \quad
k=0,1,\ldots,2j+1\;, \ea \right.
\]
and moreover we set $\tilde u_{2j+2}:=0$. The evolution of
$\Psi(\bar z(t))$ obeys the law below (a pictorial representation
of the law is given by the directed graph in Figure
\ref{fig.multivalued2}), where the symbol $\wedge$ denotes the logical conjunction `and':\\
\be\label{Psi.law}\ba{lcll} |\Psi(\bar z(t))|=\tilde u_k &\wedge &
|\bar z(t)|=\dst\frac{\tilde u_k}{1+\delta}& \Rightarrow
|\Psi(\bar z(t^+))|=\tilde u_{k+1},\, \mbox{for}\;k=0,1,\ldots,
2j+1\\[3mm]
|\Psi(\bar z(t))|=\tilde u_k & \wedge & |\bar
z(t)|=\dst\frac{\tilde u_k}{1-\delta}& \Rightarrow |\Psi(\bar
z(t^+))|=\tilde u_{k-1},\, \mbox{for}\;k=1,2,\ldots,
2j+1\\[4mm]
|\Psi(\bar z(t))|=\tilde u_{k} &\wedge & |\bar z(t)|=\tilde
u_{k-1} &\Rightarrow |\Psi(\bar z(t^+))|=\tilde u_{k-1},\,
\mbox{for}\;k=2j+2. \ea\ee If none of the conditions on the
left-hand side of the implications above is satisfied, then
$\Psi(\bar z(t^+))=\Psi(\bar z(t))$.
Observe that (\ref{Psi.law}) takes into account both the positive
and the negative values of $\Psi(\bar z(t))$. In fact, since
$\Psi(u)u\ge 0$ for all $u$, if $\Psi(\bar z(t))>0$ (respectively,
$\Psi(\bar z(t))<0$) so is $\bar z(t)$ and $\Psi(\bar z(t^+))$.
Hence, (\ref{Psi.law}) is in good accordance with Figure
\ref{fig.multivalued2}.
\\
We now specify the solution we adopt for the system \be\label{cl}
\dot x(t) = f(x(t)) + g(x(t)) \Psi(\bar z(t)) \ee with $t\in
[0,T)$. Set $t_0=0$, let $\Psi(\bar z(t_0))$ be as in
(\ref{ppsi0}), compute $\Psi(\bar z(t_0^+))$ according to
(\ref{Psi.law}) above, and consider the solution $x(t)$ of
\be\label{cl0} \dot x(t) = f(x(t))+g(x(t)) \Psi(\bar z(t_0^+)) \ee
starting from the initial condition $x_0 = \varphi(0)$, on the
interval $[t_0,t_1]$, where $t_1$ is a time at which $\bar z(t)$
satisfies one of the conditions which force $\Psi(\bar z(t))$ to
take a new value, provided that the solution of (\ref{cl0}) can be
extended up to $t_1$. By definition, $\Psi(\bar z(t))=\Psi(\bar
z(t_0^+))$ for all $t\in [t_0,t_1]$, and on $[t_0,t_1]$, $x(t)$ is
equivalently the solution of (\ref{cl}). Then, set $x_1 = x(t_1)$,
compute $\Psi(\bar z(t_1^+))$, and consider the solution of
\be\label{cl1} \dot x(t) = f(x(t))+g(x(t)) \Psi(\bar z(t_1^+)) \ee
starting from $x_1$, and defined on $[t_1,t_2]$, where $t_2$ is a
time at which a new transition occurs. Iterating this argument,
one finds a sequence $t_0, t_1, \ldots, t_k, t_{k+1}$ (for some
integer $k\ge 0$, and where we have conventionally set
$t_{k+1}=T$) of {\em switching} times, and the solution $x(t)$ of
(\ref{cl}) on $[0,T)$ is a $\overline C^1$ function of time such
that, for each $i=0,1,\ldots, k$, for all $t\in [t_i, t_{i+1})$,
it satisfies
\[
\dot x(t) = f(x(t))+g(x(t)) \Psi(\bar z(t_i^+))\;.
\]

\begin{remark}\label{rm1}
We now explain why chattering is avoided  in the interval $[0,T)$
thanks to the introduction of additional levels in the quantizer
(see also \cite{hayakawa.et.al.acc06}). In the proof of the main
result below it is shown that this property is true for all the
times. As a matter of fact, by the definition of (\ref{ppsi}),
each time $\Psi(\bar z(t))$ makes a transition from one value to
another, some (dwell) time will elapse before a new transition can
occur\footnote{For some classes of nonlinear systems, it is
possible to estimate  a lower bound on such a dwell time
(\cite{depersis.scl09}). This is particularly important in the
case in which the quantized controller is implemented over a
network, since it gives indications on the data-rate needed to
transmit the quantized information.}. This can be illustrated with
the help of Fig.\ \ref{fig.multivalued}, where $u$ is replaced by
$\bar z(t)$. Suppose that, at time $t$, $\Psi(\bar z(t))=u_0$ and
$\bar z(t)$ hits the point $\frac{u_0}{1+\delta}$ . Then
$\Psi(\bar z(t))$ takes the new value $\frac{u_0}{1+\delta}$ (see
Fig.\ \ref{fig.multivalued}). After the switching, the function
$\bar z(t)$ can increase and eventually hits the point
$u_0(1-\delta^2)^{-1}$, or decrease and eventually hits the point
$u_0(1+\delta)^{-2}$
(if it hits none of the two points then this means that $\bar
z(t)$ remains in the interval
$(u_0(1+\delta)^{-2},u_0(1-\delta^2)^{-1})$ for the entire
interval $[t,T)$, and no switching occurs in this interval).
In either case, before a new transition takes place, some time
will elapse, because the function $\bar z(t)$ must cover an
interval of finite length with finite speed. In fact,  for a given
initial condition   $\varphi\in C^1([-2\tau,0], \R^n)$, with
$||\varphi||_c\le R$ and $R>0$, the time derivative of $\bar z(t)
= z(\varphi(t-\tau))$ is continuous and bounded on $[0,T)$, and in
particular:
\[
|\frac{d\bar z(t)}{dt}|\le \max_{|x|\le R} | \frac{\partial
z(x)}{\partial x} | \cdot \max_{t\in[-2\tau, -\tau]} | \frac{d
\varphi(t)}{dt} |\;.
\]
If, on the other hand, we were adopting the quantizer
(\ref{ppsi.elia.mitter}), $\Psi(\bar z(t))$ would have taken the
value $u_1$ rather than $\frac{u_0}{1+\delta}$. Immediately after
the switching, it could happen that  $\bar z(t)$ cannot decrease,
thus forcing a transition to the previous value, which would in
turn trigger a new transition to $u_1$, and this would continue to
happen again and again. It is precisely to avoid such fast
transitions that new quantization levels were added. This addition
can be seen as a way to add {\em hysteresis} to the quantized
system, and we will refer to (\ref{ppsi}) as a quantizer with
hysteresis.
\end{remark}
For the analysis to follow, the following observation is
important. For each $t\in [0,T)$,  such that $t\in [t_i,
t_{i+1})$, $i=0,1,\ldots,k$, if $|\bar z(t)|< u_0(1-\delta)^{-1}$,
then the solution $x(t)$ of (\ref{cl}) satisfies the differential
inclusion \be\label{temfb} \dot x(t)\in f(x(t))+g(x(t))
K(\Psi(\bar z(t))) \; , \ee where $K(\Psi(u))$, with $u= \bar
z(t)$, is such that \be\label{set.k.finite} K(\Psi(u)) \subseteq
\left\{\ba{l} \{v\in \R\,:\, v = (1+\lambda\delta)u\,,\,
\lambda\in [-1,1]\} \;,
\\ \quad \quad \quad \quad \quad \quad \quad \quad \quad
(1+\delta)^{-1}u_j<|u|\le (1-\delta)^{-1}u_0\\[2mm]
\{v\in \R\,:\, v = \lambda(1+\delta)u\,,\, \lambda\in [0,1]\} \;,
 |u|\le (1+\delta)^{-1}u_j\;.
\ea\right. \ee This is easily verified bearing in mind that, by
the definition (\ref{ppsi}) of the map $\Psi(u)$, $\Psi(u)\in
K(\Psi(u))$ for all $|u|< u_0(1-\delta)^{-1}$.

\subsection{Problem formulation}

Since the control action is zero in the vicinity of the origin due
to the dead-zone of the quantizer ($\Psi(u)=0$ for $|u|\le
u_j(1+\delta)^{-1}$), asymptotic stability of the origin of
(\ref{system.ol}) is not possible to achieve (except in
exceptional cases without interest). We are rather interested in
the following property:

\begin{definition}\label{def.1}
The system \be\label{system.f.ol} \dot x(t) = f(x(t)) + g(x(t))
v(t-\tau)\;, \ee with $\tau \geq 0$ is semi-globally practically
stabilizable by quantized feedback if for any $\varepsilon < R <
0$ there exist a law $z(x)$, a real number $u_0 > 0$ and an
integer $j \geq 1$ such that the  solution of \be \label{systemf}
\dot x(t) = f(x(t)) + g(x(t)) \Psi (z(x(t-\tau))) \; , \ee
starting from ${\cal R} = \{\varphi\in C^1([-2\tau,0], \R^n):
||\varphi||_c\le R\}$ enters $B_\varepsilon$, the closed ball of
radius $\varepsilon$, at some finite time $t_s \geq 0$, and
remains in that set for all $t \ge t_s$.
\end{definition}

In the remaining sections, we propose a solution to the problem
formulated above.

\begin{remark}
The difficulty to achieve asymptotic stability can be seen by
rewriting the system (\ref{systemf}) in the form of a nominal
stable system affected by a perturbation, namely
\[
\dot x(t) = f(x(t))+ g(x(t))z(x(t)) + g(x(t))[ \Psi
(z(x(t-\tau)))- z(x(t))] \;,
\]
and neglecting the effect of the delay (the presence of the delay
worsens the situation). Consider the situation in which $z(x(t))$,
the quantity which undergoes quantization, is close to zero,
namely $|z(x(t))|< (1+\delta)^{-1} u_j$. Bearing in mind
(\ref{set.k.finite}), the perturbation $|\Psi (z(x(t)))- z(x(t))|$
is bounded from above by $|\lambda(\delta+1)-1|\,|z(x(t))|$, with
$\lambda\in[0,1]$ (the argument will be made clearer later on).
Even in the easy case in which the system is exponentially stable,
asymptotic stability cannot be proven unless the perturbation (in
this case $[\lambda(\delta+1)-1]z(x(t))$) is bounded by a linear
term $\gamma |x(t)|$ {\em and} $\gamma$ is sufficiently small (see
e.g.\ \cite{khalil.book.96}, Section 5.1), conditions which are
not met in our scenario. For the majority of the systems, these
conditions are not satisfied either  and other notions of
stability have been introduced. A notion of stability for
solutions of systems affected by non-vanishing perturbations is
that of {\em uniform ultimate boundedness} (\cite{hahn.book.67},
\cite{krasowskii.book.63}) which has found wide application in the
area of robust control (see e.g. \cite{khalil.book.96}). The
notion of semi-global practical stabilizability we consider in our
paper has been extensively investigated for problems of robust
stabilization of nonlinear systems (see e.g.\
\cite{teel.praly.sicon96}, \cite{isidori.book.2.99}, Chapter 12,
and references therein). The same notion of stability has been
already studied for quantized time-delay systems as well
(\cite{liberzon.tac06}, \cite{teel.tac98}).
\end{remark}

\section{Standing assumptions and main result}
\label{mainre}

\subsection{Basic assumptions}\label{3.1}

The result to be derived below for the system (\ref{system.f.ol})
holds
under the following standing assumptions.\\
\noindent {\bf (A1)} There exist a continuously differentiable
positive definite and proper Lyapunov function $V(x)$, two class
${\cal K}_\infty$ functions $\kappa_1, \kappa_2$, a positive
definite continuous function $W(x)$ and a continuously
differentiable real-valued function  $z(x)$, which is zero at the
origin, with $W(x)$ and $z(x)$ both depending on $\delta$, such
that, for all $x \in \R^n$, \be\label{nom.stability} \ba{cl}
\kappa_1(|x|)\le V(x)\le \kappa_2(|x|) \ , &
\\
\dst\frac{\partial V}{\partial x}[f(x) + g(x)(1+p)z(x)] \le - W
(x)\;, & p\in [-\delta,\delta]\;. \ea \ee

\begin{remark}
It would be slightly more correct to denote $W(x)$ and $z(x)$ by,
respectively, $W_\delta(x)$ and $z_\delta(x)$, since due to the
presence of the uncertainty in the input channel, both these
functions are going to depend on the size $\delta$ of the
uncertainty (see Subsection \ref{subs.a1} below). However, to ease
the notation, we decided not to make the dependence on $\delta$
explicit.
\end{remark}

\begin{remark}
The uncertainty in the input channel is modeled through the
parameter $p$, whose range depends on the quantization density
through $\delta$. Such uncertainty takes into account the effect
due to quantization, as it should be evident from
(\ref{set.k.finite}). Assumption (A1) amounts to require the
system $\dot x(t)=f(x(t))+g(x(t))u(t)$, with {\em no} delay,  to
be stabilizable in the presence of quantization. The design of a
stabilizing quantized feedback is carried out e.g.\ in
\cite{ceragioli.depersis.scl07} (see also Subsection \ref{subs.a1}
below).
\end{remark}

The next two assumptions require the system to be robust with
respect to delays. In particular they are needed to guarantee that
no finite-escape time phenomenon will occur, and that the solution
stays bounded for all the times. These conditions also appear in
\cite{mazenc.bliman.tac06} (where no quantization was present),
although in a slightly different form. The difference is due to
the fact that the quantization effect adds up to the delay effect,
and in the conditions below also the quantization parameter
$\delta$ plays a role. More comments on these two assumptions are
postponed to Subsection \ref{subs.a2.a3}.

\medskip

\noindent {\bf (A2)} Let $\Omega$ be a positive real number which
satisfies $\Omega\ge 16 \tau$. For all $x\in \R^n$, for all $\xi
\in \overline C^1([0,2\tau],\R^n)$, for all $\lambda_1 \in [-1,1]$
and for all $\lambda_2 \in \overline C^0([0,\tau],\R)$ such that
$\lambda_2(m) \in [1-\delta,1+\delta]$ for all $m \in [0,\tau]$,
the inequality \be \label{in.h2} - \dst\frac{1}{4} W (x) -
T(x,\xi,\lambda_1,\lambda_2)-\dst\frac{1}{\Omega}\dst\int_0^{2\tau}
W (\xi(\ell)) d\ell \le 0\;, \ee with
\[\ba{rcl}
T(x,\xi,\lambda_1,\lambda_2) & = & L_gV(x)(1+\lambda_1\delta)
\dst\int_{\tau}^{2\tau} H(\xi(\ell), \xi(\ell-\tau),
\lambda_2(\ell-\tau)) d\ell \; ,
\\
H(a,b,c) & = & L_f z(a) + L_g z(a) c z(b)\;, \ea\] holds.

\noindent {\bf (A3)} There exists a nondecreasing function
$\kappa_{3}(\cdot)$ of class $C^1$ such that for all $x \in \R^n$,
for all $L \geq 0$ and for all $\lambda \in [-1,1]$, the
inequality \be\label{ineq.h3'} - \dst\frac{1}{2}W (x) +
\displaystyle\sup_{|a| \leq L}\left\{L_g V(x) (1+\lambda)[z(a) -
z(x)]\right\} \le \kappa_{3}(L) [V(x) + 1] \ee holds. Let
$\kappa_4(L) = 2\kappa_{3}(L)$.

\medskip

\noindent The two subsections below provide comments to help the
readers to understand the role played by each assumption in the
solution of the problem. However, the reader who is interested in
getting to the statement of the main result immediately, can skip
the next two subsections and go directly to Subsection
\ref{subs.main.result}.

\subsection{Comments on the Assumption (A1)}\label{subs.a1}

A number of ways to have Assumption (A1) fulfilled are discussed
below.

\begin{itemize}
\item {\em Lyapunov Redesign}. Suppose that,
for the system (\ref{system.f.ol}), are known a function $V$ of
class $C^2$, and a function $\zeta(x)$ of class $C^1$ such that,
instead of (\ref{nom.stability}), only the weaker condition
\be\label{nom.stability2} L_f V(x) + L_g V(x) \zeta(x) = - \tilde
W(x) \;, \ee with $\tilde W(x)$ a continuous positive definite
function, is satisfied. Introduce the control law \be\label{zeta}
z(x) = \zeta(x) - \alpha(x) L_g V(x) \; , \ee with $\alpha(x)$ a
positive function to be chosen later. Then we have
\[\ba{l}
\dst\frac{\partial V}{\partial x}[f(x)+g(x)(1+p) z(x)]\\[2mm]
 =  \dst\frac{\partial V}{\partial x}[f(x) + g(x) \zeta(x)]
- \alpha(x)\left|L_g V(x)\right|^2+ p L_g V(x)
\left[\zeta(x) -\alpha(x) L_g V(x)\right]\\[2mm]
 \le  -W(x) - \alpha(x)(1+p)\left|L_g V(x)\right|^2 + p L_g V(x)\zeta(x)\\
 \le  -W(x) - \alpha(x)(1-\delta)\left|L_g V(x)\right|^2+
\delta\left|L_g V(x)\right||\zeta(x)| \;. \ea\] A simple
completion-of-the-squares argument shows that
\[
\dst\frac{\partial V}{\partial x}[f(x)+g(x)(1+p)z(x)] \le
-\dst\frac{3}{4} \tilde W(x)\;,
\]
provided that \be\label{alpha} \alpha(x)\ge
\dst\frac{\delta^2}{1-\delta}\dst\frac{|\zeta(x)|^2}{\tilde
W(x)}\;. \ee Hence, the control law (\ref{zeta}), with $\alpha(x)$
defined above and such that \linebreak $\lim_{x\to 0}
\alpha(x)L_gV(x)=0$, guarantees the fulfillment of Assumption (A1)
with $W (x)=3\tilde W(x)/4$.

\item {\em Sontag's universal stabilizer \cite{sontag.scl89}.} Consider
the system \be \label{sysap} \dot x = f(x) + g(x) [1 + p] u \; ,
\ee with $x \in \R^n$, $u \in \R$, $p \in [-\delta, \delta]$,
$\delta \in [0, 1)$. Let us assume that a control Lyapunov
function $V(x)$ is known for the system (\ref{sysap}) with $p =
0$, and set \be \label{z1} \dot V(x) = a(x) + [1 + p] b(x) u \; ,
\ee with \be \label{z2} a(x) = L_f V(x) \; , \; b(x) = L_g V(x) \;
. \ee Since $V$ is a control Lyapunov function for (\ref{sysap})
with $p = 0$, $b(x) = 0$ implies $a(x) < 0$ when $x \neq 0$. Next,
consider the control given by Sontag's formula: \be
\label{z3}\ba{l} u(x) = K \frac{- a(x) - \sqrt{a(x)^2 +
b(x)^4}}{b(x)} \; \mbox{when} \; b(x) \neq 0
\; ,\\
u(x) = 0 \; \mbox{when} \; b(x) = 0\;, \ea \ee and where $K$ is a
positive real number to be selected later. Then, when $b(x) \neq
0$, the derivative of $V$ along the trajectories of (\ref{sysap})
in closed-loop with $u(x)$ defined in (\ref{z3}) satisfies
\begin{equation}
\begin{array}{rcl}
\label{z4} \dot V(x) & = & a(x) + [1 + p] b(x) K \dst\frac{- a(x)
- \sqrt{a(x)^2 + b(x)^4}}{b(x)}
\\[2mm]
& = & a(x) - [1 + p] K a(x) - [1 + p] K \sqrt{a(x)^2 + b(x)^4}
\\[2mm]
& = & [1 - (1 + p) K]a(x) - [1 + p] K \sqrt{a(x)^2 + b(x)^4}\;.
\end{array}
\end{equation}
We choose $K = \frac{2}{1 - \delta} > 0$. Then, when $a(x) \geq
0$, we have
\begin{equation}
\begin{array}{rcl}
\label{z5} \dot V(x) & \le & - a(x) - [1 + p] K \sqrt{a(x)^2 +
b(x)^4}\;,
\end{array}
\end{equation}
and, when $a(x) < 0$,
\begin{equation}
\begin{array}{rcl}
\label{z6} \dot V(x) & = & a(x) - [1 + p] K (a(x) + \sqrt{a(x)^2 +
b(x)^4}) < 0 \; .
\end{array}
\end{equation}
When $b(x) = 0$, then
\begin{equation}
\begin{array}{rcl}
\label{z7} \dot V(x) & = & a(x) < 0 \; \mbox{if} \; x \neq 0\;.
\end{array}
\end{equation}
Under the {\em small control property} (\cite{sontag.scl89}) one
can prove that the control law introduced above is smooth
everywhere except at the origin where it may be only continuous.
However, in many cases, the control law turns out to be also
continuously differentiable at the origin, and then a continuously
differentiable function $z(x)$ which guarantees the inequality
(\ref{nom.stability}) is obtained.

\item {\em Lyapunov stable systems.}
Consider again system (\ref{sysap}), and assume that a Lyapunov
function $V(x)$ such that $a(x) \leq 0$, and $b(x) \neq 0$ when $x
\neq 0$ and $a(x) = 0$, is known. Then, selecting \be \label{s1} u
= - \xi(x) b(x) \;, \ee where $\xi$ is any $C^1$ positive
function, we obtain, for all $x \in \mathbb{R}^n$ \be \label{z11}
\dot V(x) \leq a(x) - [1 - \delta] \xi(x) b(x)^2 \ee and the
function $a(x) - [1 - \delta] \xi(x) b(x)^2$ is negative definite.
Inequality (\ref{nom.stability}) then holds with $z(x)= - \xi(x)
b(x)$ and $W (x)=-a(x) + [1 - \delta] \xi(x) b(x)^2$.

\medskip

\item {\em Dissipation inequality \cite{hayakawa.et.al.acc06},
\cite{ceragioli.depersis.scl07}}. Consider the system
(\ref{system.f.ol}). Suppose that  a Lyapunov function $V(x)$ is
known such that for all $x \in \mathbb{R}^n$
\[
L_f V(x) -\dst\frac{1}{4}(1-\delta^2)\left(L_g V(x)\right)^2
 \le - \tilde W(x)\;.
\]
Then, for any $p\in [-\delta, \delta]$ it is also true that
\[
 L_f V(x)
-\dst\frac{1}{4}(1-p^2) \left(L_g V(x)\right)^2  \le -\tilde
W(x)\;.
\]
Define now
\[
z(x) = - \dst\frac{1}{2} L_g V(x)
\]
and observe that the inequality above rewrites as
\[
L_f V(x) + \dst\frac{1}{4}p^2 \left(L_g V(x)\right)^2 + z(x) L_g
V(x) + z(x)^2 \le - \tilde W(x)\;,
\]
or, equivalently, \be\label{hji} \dst\frac{\partial V}{\partial x}
(f(x)+g(x) z(x))+\dst\frac{1}{4}p^2 \left(L_g V(x)\right)^2 +
z(x)^2\le -\tilde W(x) \; . \ee We remark incidentally
(\cite{ceragioli.depersis.scl07}) that the latter inequality
implies the existence of a control $u= z(x)$ which renders the
system
\[
\left\{ \ba{rcl} \dot x & = & f(x) + g(x)u + g(x)w \ ,
\\
z & = & u \ , \ea \right.
\]
{\em strictly dissipative} with respect to the supply rate
$q(w,z) = - z^2 + p^{-2} w^2$.\\
Observe now that
\[
p z(x) L_g V(x) \le \dst\frac{1}{4}p^2 \left(L_g V(x)\right)^2 +
z(x)^2
\]
and therefore (\ref{hji}) implies that
\[
\dst\frac{\partial V}{\partial x} (f(x)+g(x)z(x)) + p  z(x) L_g
V(x) \le - \tilde W(x)\;,
\]
that is (\ref{nom.stability}) with $W (x)=\tilde W(x)$.
\end{itemize}

\subsection{Comments on the Assumptions (A2) and (A3)}\label{subs.a2.a3}
The two Assumptions (A2) and (A3) describe, in terms of the
Lyapunov function $V$, how robust with respect to delays in the
input channel the system $\dot x(t)=f(x(t))+g(x(t))u(t)$ should be
in order to find a stabilizing feedback despite the delay. The
role of these assumptions for systems with no quantization was
already investigated in \cite{mazenc.bliman.tac06}. To better
assess such a role, let us neglect the effect due to the
quantization, and  let us set $\delta=0$. Then, in Assumption
(A1), $p=0$ and (\ref{nom.stability}) becomes a standard
stabilizability assumption. The inequality (\ref{in.h2}) in (A2)
becomes \be \label{in.h2.2} - \dst\frac{1}{4} W (x) -
T(x,\xi)-\dst\frac{1}{\Omega}\dst\int_0^{2\tau} W (\xi(\ell))
d\ell \le 0\;, \ee with
\[\ba{rcl}
T(x,\xi) & = & L_gV(x) \dst\int_{\tau}^{2\tau} H(\xi(\ell),
\xi(\ell-\tau)) d\ell \; ,
\\
H(a,b) & = & L_f z(a) + L_g z(a) z(b)\;. \ea\] Similarly, in
(\ref{ineq.h3'}), $\lambda=0$, and the inequality implies that for
all $\xi\in {\cal C}^1([0,2\tau], \R^n)$, there exists a positive
constant $\kappa_\xi$ such that, for all $x \in \R^n$, for all
$t\in [0,2\tau]$, \be\label{ineq.h3''} - \dst\frac{1}{2}W (x) +
L_g V(x)[z(\xi(t)) - z(x)] \le \kappa_\xi [V(x) + 1]\;. \ee The
conditions (\ref{in.h2.2}), (\ref{ineq.h3''}) coincide with those
found in \cite{mazenc.bliman.tac06} to prove that the origin of
\be\label{cl.no.quantization} \dot
x(t)=f(x(t))+g(x(t))z(x(t-\tau)) \ee is uniformly globally
asymptotically stable. Compared with \cite{mazenc.bliman.tac06},
the stronger conditions we have in this paper are due to the fact
that both quantization and
delay affect the system.\\
In the case no quantization is present, the role of
(\ref{in.h2.2}), (\ref{ineq.h3''}) to guarantee stability of
time-delay systems is easier to describe (see
\cite{mazenc.bliman.tac06} for details). The condition
(\ref{ineq.h3''}), for instance, guarantees that no finite-time
escape of the solution occurs. As a matter of fact, the time
derivative of $V$ computed along the solutions of
(\ref{cl.no.quantization}) obeys the equations
\[\ba{rcl}
\dot V(x(t)) &=& L_f V(x(t))+ L_g V(x(t))z(x(t-\tau))\\
&=& L_f V(x(t))+ L_g V(x(t))z(x(t))+L_g V(x(t))[z(x(t-\tau))-z(x(t))]\\
&\le& -W(x(t))+ L_g V(x(t))[z(x(t-\tau))-z(x(t))]\;. \ea\] As $t$
ranges in the interval $[0, \tau]$,  $x(t-\tau)$ can be viewed as
a function $\xi\in {\cal C}^1([-\tau,0], \R^n)$, and bearing in
mind (\ref{ineq.h3''}), we have $\dot V(x(t))\le K_\xi
(V(x(t))+1)$. From this we infer that no finite escape-time can
exist on $[t_0, t_0+\tau]$. Iterating the argument, one can prove
that the solution is defined for all $t$.\\
The condition (\ref{in.h2.2}) guarantees that a suitable
Lyapunov-Krasowskii functional is strictly decreasing along the
solutions of the closed-loop system (again, the interested reader
is referred to \cite{mazenc.bliman.tac06} for more details). The
purpose of the rest of the paper is to show how, taking advantage
of Assumptions (A1)-(A3), the arguments of
\cite{mazenc.bliman.tac06} can be modified to take into account
the additional constraints due to the presence of the
quantizer.\\
We stress that the conditions  (\ref{in.h2.2}), (\ref{ineq.h3''})
require the system to be robust with respect to (quantization and)
delays and are essential to design stabilizing control laws for
nonlinear (quantized) time-delay systems. Analogous conditions are
found in other contributions on the topic. In \cite{teel.tac98},
using an approach based on Razumikhin-like theorems, uniform
asymptotic stability with restriction $\Delta$ on the norm
$||\varphi||_c$ of the initial condition and with offset
$\varepsilon$ (a notion of stability very similar to what we have
in Definition \ref{def.1}) is proven.
%
To be more precise, suppose the system is stabilizable, that is
(A1) holds (with $p=\delta=0$). Also suppose for the sake of
simplicity that $W(x)$ is replaced by the class-${\cal K}_\infty$
function $\alpha_3(|x|)$. Then it is possible to design a smooth
invertible function $G(x)$ and a class-${\cal K}_\infty$ function
$\gamma_\theta$ such that \be\label{one} |x(t)|\ge
\gamma_\theta(|\theta(t)|)\; \Rightarrow\; \dot V(t)\le
-\dst\frac{1}{2}\alpha_3(|x(t)|)\;, \ee with
\[\ba{rcl}
\theta(t)&=&-G^{-1}(x(t))[z(x(t))-z(x(t-\tau))]\\
&=&-G^{-1}(x(t)) \dst\int_{t-\tau}^t\left.\dst\frac{\partial
z(x)}{\partial x}\right|_{x=x(s)} [f(x(s))+g(x(s))z(x(s-\tau))]ds.
\ea\] Further, one can find a class-${\cal K}$ function $\gamma_1$
such that \be\label{two} |\theta(t)|\le \tau
\gamma_1(\sup_{t-2\tau\le s\le t} |x(s)|)\;. \ee Hence, combining
(\ref{one}) and (\ref{two}), one obtains:
\[
|x(t)|\ge \gamma_\theta(\tau \gamma_1(\sup_{t-2\tau\le s\le t}
|x(s)|))\; \Rightarrow\; \dot V(t)\le
-\dst\frac{1}{2}\alpha_3(|x(t)|)\;.
\]
Under the small-gain condition \be\label{small.gain}
\kappa_1^{-1}\circ \kappa_2\circ
\gamma_\theta(\tau\gamma_1(s))<s\;,\quad \mbox{for all}\;
\varepsilon<s<\Delta\;, \ee the inequality above shows that the
zero solution of the system
\[
\dot x(t)=f(x(t))+g(x(t))\zeta(x(t-\tau))
\]
is uniformly asymptotically stable with restriction $\Delta$ on
the norm $||\varphi||_c$ of the initial condition, and with offset
$\delta$ \footnote{In the terminology of  [21], the offset is the
size of the set where the state converges at some finite time and
stays there from that time on -- in our paper such a
parameter is denoted by $\varepsilon$.}.\\
The condition (\ref{small.gain}) represents an alternative way to
express robustness of the system with respect to delays to infer
stability results using Razumikhin-like theorems.

\subsection{Main result}\label{subs.main.result}
We are ready to state the main result of our work. As already made
clear in the problem formulation (Definition \ref{def.1}), the two
main design parameters are the range $u_0$ and the number $j$ of
levels of the quantizer. Intuitively, to design  $u_0$ we need to
quantify the ``overshoot" of the state variable and we expect this
to depend on the size of the initial condition. Regarding the
number of quantization levels $j$, it is not hard to figure out
that in general the closer one wants to confine the state to the
origin (i.e.\ the smaller $\varepsilon$ is in Definition
\ref{def.1}), the larger the number of quantization levels must
be. On the other hand, having fixed the width of the quantizer,
the  number of the quantization levels will increase with the
range $u_0$ and in turn with $R$. Such a dependence is made clear
in the statement below. The proof is constructive and provides the
explicit expressions for  $u_0$ and $j$.

\begin{proposition}
\label{p2} Let us assume that the system (\ref{system.f.ol})
satisfies Assumptions (A1) to (A3). Then the origin of
(\ref{system.f.ol}) is semi-globally practically stabilizable by
quantized feedback. Namely, there exist a positive, continuous and
non-decreasing function $u_0(\cdot):\R_{\ge 0}\to \R_{> 0}$, and
a positive continuous function $j(\cdot,\cdot):\R_{\ge 0}^2\to
\R_{> 0}$ such that, for any $R > \varepsilon > 0$, if $u_0\ge
u_0(R)$,  $j\ge j(\varepsilon,R)$ and $z$ is the feedback provided
by Assumption (A1) satisfying (\ref{nom.stability}), then the
solution of (\ref{systemf}) starting from ${\cal R}=\{\varphi\in
C^1([-2\tau,0], \R^n): ||\varphi||_c\le R\}$ enters
$B_\varepsilon$, the closed ball of radius $\varepsilon$, at some
finite time $t_s \geq 0$, and remains in that set for all $t \ge
t_s$.
\end{proposition}

\noindent The proof of the result is postponed to the next
section. Before ending Section \ref{mainre}, we discuss two
examples in which the proposition above is applied.

\subsection{Example 1}

We illustrate Proposition \ref{p2} by showing how it applies when
the functions $f$ and $g$ in (\ref{system.f.ol}) are linear. Thus,
we consider the system \be\label{lin.sys} \dot x(t) = Ax(t) +
B\Psi (u(t-\tau))\;, \ee where $A \in \mathbb{R}^{n \times n}$ and
$B \in \mathbb{R}^{n \times 1}$ are constant matrices. We assume
that the pair $(A,B)$ be stabilizable. Then there exist a positive
definite symmetric matrix $Q \in \mathbb{R}^{n \times n}$ and a
matrix $\tilde{K} \in \mathbb{R}^{1 \times n}$ such that
\[
(A + B\tilde K)^\top Q + Q(A + B\tilde K) = - I
\]
where $I \in \mathbb{R}^{n \times n}$ denotes the identity matrix.
Then, in view of the Lyapunov redesign we have proposed to
determine a control law such that Assumption (A1) is verified, one
can verify that the matrix
\[
K=\tilde K-2\alpha Q B\;,\;{\rm with}\; \alpha\ge
\dst\frac{\delta^2}{1-\delta}|K|^2\;,
\]
is such that, for all $x \in \mathbb{R}^n$, \be\label{a1.rld}
2x^\top Q [Ax + B(1+p)K x]\le -c x^Tx\;,\quad p\in
[-\delta,\delta] \ee with $c = 3/4$. Therefore Assumption (A1) is
satisfied with $V(x)=x^\top Q x$, $z(x) = Kx$, and $W (x) = c
x^\top x$.
Hence in what follows we let $Q,K,c$ be such that (\ref{a1.rld}) holds. \\
We turn now to Assumption (A2). We have
\[\ba{l}
-T(x,\xi,\lambda_1,\lambda_2)  =  - 2(1+\lambda_1\delta)
x^TQB\dst\int_\tau^{2\tau} K
\{A\xi(l)- B \lambda_2(\ell-\tau) K \xi(\ell-\tau)\}d\ell\\
\le 2(1+\delta) |x|\,|QB|\dst\int_\tau^{2\tau} |KA|
|\xi(\ell)|d\ell
+ 2(1+\delta)^2 |x|\,|QB|\dst\int_\tau^{2\tau} |KBK|\cdot \\
 \hspace{10cm} \cdot|\xi(\ell-\tau)|d\ell\\
\le 2(1+\delta) |x|\,|QB|\dst\int_\tau^{2\tau} [G_1 |\xi(\ell)| +
G_2 |\xi(\ell-\tau)|] d\ell \ea\] with $G_1 = |KA|$, $G_2 =
(1+\delta)|KBK|$. In view of the bounds on
$T(x,\xi,\lambda_1,\lambda_2)$, Assumption (A2) is verified if
\be\label{i2} \ba{l} -\dst\frac{c}{4} |x|^2+ 2(1+\delta)
|x|\,|QB|\dst\int_\tau^{2\tau} [G_1 |\xi(\ell)|
+ G_2 |\xi(\ell-\tau)|]d\ell\\[2mm]\hspace{6.5cm}
-\dst\frac{c}{\Omega} \dst\int_0^{2\tau} |\xi(\ell)|^2 d\ell\le 0
\ea\ee with $\Omega = 16\tau$. We easily deduce that (\ref{i2}) is
satisfied if \be\label{i2.2} \ba{l} -\dst\frac{c}{4}|x|^2 +
2G_3|x|\dst\int_0^{2\tau} |\xi(\ell)|d\ell -\dst\frac{c}{16\tau}
\dst\int_0^{2\tau} |\xi(\ell)|^2 d\ell\le 0 \ea\ee with $G_3 =
2|QB|(1+\delta)\max\{G_1, G_2\}$. By Young's inequality applied to
the second term, we deduce that (\ref{i2.2}) is satisfied if there
exists $\varepsilon > 0$ such that \be\label{it} \ba{l}
\left(-\dst\frac{c}{4} + \dst\frac{G^2_3}{\varepsilon}
\right)|x|^2 +\left(\varepsilon\cdot 2\tau
-\dst\frac{c}{16\tau}\right) \dst\int_0^{2\tau} |\xi(\ell)|^2
d\ell \le 0\;. \ea\ee The inequality holds if $\varepsilon =
4G_3^2/c$ and \be\label{tau.bound}\ba{rcl} \tau &\le &
\dst\frac{c}{8\sqrt{2} G_3}\le \dst\frac{c}{16(1+\delta)
\sqrt{2}|QB|\max\left\{|KA|, (1+\delta)|KBK| \right\}}\;. \ea\ee
Finally we consider Assumption (A3). The left-hand side of
(\ref{ineq.h3'}) becomes
\[\ba{rcl}
&&-\dst\frac{1}{2}W (x)+ \displaystyle\sup_{|a| \leq L} \{ L_g
V(x) (1+\lambda)[z(a) - z(x)] \}
\\[2mm]
& = & -\dst\frac{c}{2}x^T x +
 \displaystyle\sup_{|a| \leq L}
\{ 2x^T Q B(1+\lambda) K [a-x] \}
\\
& \le & 4 |Q B K||x|[|x| + L]\\
& \le & 4 |Q B K|[(L+1)|x|^2 + L]\\
& \le & 4 |Q B K|(L+1)[\lambda_{min}^{-1}(Q)V(x) + 1]\; . \ea\] We
deduce that one can find a constant $\Gamma = 4 |Q B
K|\max\{\lambda_{min}^{-1}(Q),1\}$ such that Assumption (A3) is
satisfied with $\kappa_3(\ell) = \Gamma (\ell + 1)$. Summarizing,
Assumptions (A1)-(A3) are satisfied for the system
(\ref{lin.sys}). Hence, we can conclude that Proposition \ref{p2}
applies, provided that the pair $(A,B)$ is stabilizable, and the
delay $\tau$ satisfies (\ref{tau.bound}).

\subsection{Example 2}

\begin{figure}
\begin{center}
\scalebox{0.75}{\includegraphics{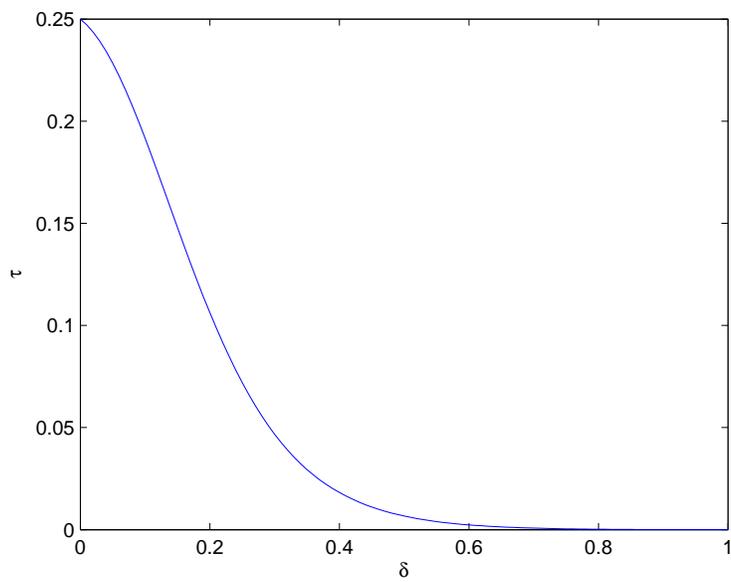}}
\caption{\label{fig.regionstability} The region under the graph
represents the set of pairs $(\delta,\tilde\tau)$ for which the
system (\ref{example2}) is semi-globally practically stabilizable.
In the picture, $\tilde \tau$ is simply denoted as $\tau$. As
$\delta$ tends to $0$ (no quantization) the normalized maximal
allowable delay $\tilde \tau$ approaches its maximum. }
\end{center}
\end{figure}

In this section we consider the classical equations of an actuated
pendulum without friction: \be\label{example2} \ba{rcl}
\dot x_1 &=& x_2\\
\dot x_2 &=& -\sin x_1 +u\;. \ea\ee The control law
\[
u=\zeta (x)=\sin x_1-x_1-2x_2
\]
and the Lyapunov function
\[
V(x)= x^T Q x = x^T \left( \ba{cc}
\dst\frac{3}{2} &  \dst\frac{1}{2}\\[2mm]
\dst\frac{1}{2} &  \dst\frac{1}{2} \ea \right) x
\]
are such that
\[
\dst\frac{\partial V}{\partial x}(f(x)+g(x)\zeta(x)) =-|x|^2\;.
\]
Applying the Lyapunov redesign of Subsection \ref{subs.a1}, it is
straightforward to see that
\[
z(x) = \sin x_1 -(\alpha+1)x_1 -(\alpha+2)x_2\;,
\]
with $\alpha\ge 16\frac{\delta^2}{1-\delta}$, guarantees
Assumption (A1) with $W(x)=\frac{3}{4}|x|^2$.
To check Assumption (A2), observe that
\[\ba{rcl}
&&-T(x,\xi,\lambda_1,\lambda_2)\\
&\le& 2(1+\delta)|x|\dst\int_{\tau}^{2\tau}
\left[2(2+\alpha)|\xi(\ell)|+
(2+\alpha)(1+\delta)2(2+\alpha)|\xi(\ell-\tau)|\right] d\ell\\
&\le& 2(1+\delta)|x|2(2+\alpha)\dst\int_{\tau}^{2\tau}
\left[|\xi(\ell)|+ (2+\alpha)(1+\delta)|\xi(\ell-\tau)|\right]
d\ell\;. \ea\] Similarly to the previous example, one can prove
that, if
\[
\tau\le \dst\frac{3}{128\sqrt{2}(\alpha+2)^2(1+\delta)}\;,
\]
then Assumption (A2) is fulfilled. Even Assumption (A3) can be
easily verified. As a matter of fact,
\[\ba{rcl}
&&- \dst\frac{1}{2}W (x) + \displaystyle\sup_{|a| \leq
L}\left\{L_g V(x) (1+\lambda)[z(a) -
z(x)]\right\}\\
&\le & \displaystyle\sup_{|a| \leq L}\left\{2|x|
(1+\delta)2(\alpha+2)[|a|+|x|]
\right\}\\
&\le & 4 (1+\delta)(\alpha+2)|x|[|x|+L]\;. \ea
\]
As in the previous example, one can deduce that Assumption (A3) is
fulfilled with $\kappa_3(\ell)=\Gamma (\ell+1)$ and $\Gamma=4
(1+\delta)(\alpha+2)\frac{\sqrt{2}}{\sqrt{2}-1}$. The region under
the graph in Fig.\ \ref{fig.regionstability} describes the pairs
$(\delta,\tilde\tau)$, with $\tilde\tau=
\frac{1}{(\alpha+2)^2(1+\delta)}$, for which the system
(\ref{example2}) is semi-globally practically stabilizable.

\section{Proof of Proposition \protect\ref{p2}}\label{s.proof}

The proof is based on a Lyapunov-Krasowskii functional given by
the sum of the Lyapunov function $V(x)$ in Assumption (A1) and a
term which at time $t$ depends on the state $x(\cdot)$ restricted
to the interval $[t-2\tau,t]$. Hence, in order to use such a
Lyapunov-Krasowskii functional, we need to
first prove that all solutions of the closed-loop system we
consider exist for all $t \in [- 2\tau, 2\tau]$. To this purpose,
we will only make use of the Lyapunov function $V(x)$. Then we
will prove that the solutions
can be extended beyond $2\tau$, showing that the
Lyapunov-Krasowskii functional is bounded for all the time
and finally that the solutions converge in finite time to a ball
around the origin of radius $\varepsilon$.
\subsection{Existence of solutions for $t\in [-2\tau,2\tau]$}

As a first step, we need to define the  function $u_0(\cdot)$ by
which we define the range $u_0$. We have already observed that to
find such a function, we need to estimate the region where the
state is confined for all the times. We will obtain such an
estimate by steps, first estimating a bound on $|x(t)|$ on the
interval $[0, \tau]$, then a bound on the interval $[0, 2\tau]$,
and finally a bound on $[0,+\infty)$. Let us then introduce such
sequence of bounds as functions of the nonnegative real-valued
parameter $R$, the radius of the ball of initial conditions:
\be\label{und} \alpha(R) = \kappa_1^{-1}\left(e^{\kappa_4(R)
\tau}(\kappa_2(R) + 1) - 1\right) \; , \ee
\be\label{oune} \gamma(R) = \alpha(\alpha(R)) \; , \ee
\be\label{berl} \omega(R) = \kappa_1^{-1}\left(\kappa_2(\gamma(R))
+ \frac{\tau}{4} \displaystyle\sup_{|a| \leq \gamma(R)}
W(a)\right) + R \; , \ee where $\kappa_1, \kappa_2,\kappa_4$ are
the class ${\cal K}_\infty$ functions defined in Section
\ref{3.1}. Observe that the functions $\alpha, \gamma, \omega$ are
continuous and for all $R \geq 0$, the inequalities
\be\label{ozbis} \omega(R) \geq \gamma(R) \geq \alpha(R) \geq R
\ee are satisfied. It will be proven below that $|x(t)|\le
\omega(R)$ for all $t\ge -2\tau$. Define \be\label{u02} u_0(R) =
\displaystyle\sup_{|a| \leq \omega(R)}|z(a)| + 1 \ee and let
$u_0\ge u_0(R)$. Having defined $u_0$ we can proceed with the rest
of the proof.
\\
Consider the  solution $x(t)$ of (\ref{systemf}) with an initial
condition $\varphi\in {C}^1([-2 \tau,0], \R^n)$ such that
$||\varphi||_c \le R$. Let us show first that this solution is
defined over $[- 2\tau,\tau]$. To prove this, let us proceed by
contradiction. Suppose it is not defined over $[- 2\tau,\tau]$.
Observe that, since $||\varphi||_c \le R$, then $|z(\varphi(t -
\tau))| < (1 - \delta)^{-1}u_0$ for all $t\in [0, \tau]$, by
definition of $u_0(\cdot)$ (see (\ref{u02})). Hence, $\Psi
(z(\varphi(t - \tau)))$ is well-defined for all $t\in [0, \tau]$.
Next, we deduce that, necessarily there exists $T\in (0,\tau]$
such that the solution exists for all $t \in [0, T)$. Such
solution satisfies, for all $t\in [0, T)$ such that $t\in
[t_i,t_{i+1})$, $i=0,1,\ldots,k$, the differential inclusion
\be\label{fi} \dot x(t) \in f(x(t))+g(x(t)) K(\Psi (z(\varphi(t -
\tau)))) \; , \ee where $K(\Psi (z(\varphi(t - \tau))))$ denotes
the set (\ref{set.k.finite}) with $u = z(\varphi(t - \tau))$. For
all $t \in [0, T)$ such that $t\in [t_i,t_{i+1})$, $i =
0,1,\ldots,k$, we are interested in finding an upper bound for the
term \be\label{up} L_f V(x(t)) + L_g V(x(t)) v \;, \ee for any
$v\in K(\Psi (z(\varphi(t - \tau))))$. Indeed, since for any $t\in
[t_i,t_{i+1})$, $i=0,1,\ldots,k$, the derivative of $V$ along the
trajectories of the system we consider satisfies
\[\ba{rcl}
\dot V(t) & = & L_f V(x(t)) + L_g V(x(t))\Psi (z(\varphi(t -
\tau))) \ea\] then
\[
\dot V(t)= L_f V(x(t)) + L_g V(x(t)) v
\]
for some $v\in K(\Psi (z(\varphi(t - \tau))))$, and finding an
upper bound for (\ref{up}) means providing an upper bound for
$\dot V(t)$. Observe that to find an upper bound for  (\ref{up}),
it suffices to find an upper bound for \be\label{va} \ba{l} v_a :=
L_f V(x(t)) + L_g V(x(t))(1 + \lambda\delta) z(\varphi(t -
\tau))\\ \hspace{6cm} + L_g V(x(t))(\lambda - 1)z(\varphi(t -
\tau)) \varrho \ea\ee for all $t \in [0, T)$, with $\lambda$ any
number in the interval $[-1,1]$ and where $\varrho = 1$ or $0$. As
a matter of fact, since in (\ref{up}), $v\in K(\Psi (z(\varphi(t -
\tau))))$, then, by (\ref{set.k.finite}), either
\[
v=(1+\lambda\delta)z(\varphi(t - \tau))\;,\;\mbox{with}\;
\lambda\in [-1,1]
\]
(provided that $(1+\delta)^{-1}u_j<|z(\varphi(t-\tau))|\le
(1-\delta)^{-1}u_0$, in which case $\varrho=0$) or
\[
v=\lambda(1+\delta)z(\varphi(t - \tau))\;,\;\mbox{with}\;
\lambda\in [0,1]
\]
(provided that $|z(\varphi(t-\tau))|\le
(1+\delta)^{-1}u_j$, in which case $\varrho=1$).\\
Hence, for  a fixed $t$, the set of values in (\ref{up}) obtained
as $v$ ranges over  $K(\Psi (z(\varphi(t - \tau))))$ is contained
in the set of values of $v_a$ as $\lambda\in [-1,1]$ and
$\varrho\in \{0,1\}$.

Now, adding and subtracting $L_g V(x(t))(1+\lambda \delta)z(x(t))$
on the right-hand side of the equality (\ref{va}), and taking
advantage of (\ref{nom.stability}), we deduce that \be\label{ot.v}
\ba{l} v_a \leq - W (x(t)) + L_g V(x(t))
(1+\lambda\delta)[z(\varphi(t - \tau)) - z(x(t))]\\
\hspace{5.5cm} + L_g V(x(t)) (\lambda - 1) z(\varphi(t -
\tau))\varrho\; . \ea \ee Since $\delta\in (0,1)$, for any $t\in
[0,T)$ and $\lambda\in [-1,1]$, the quantity $L_g V(x(t))
(1+\lambda\delta)[z(\varphi(t - \tau)) - z(x(t))]$ belongs to  the
set
\[
\{L_g V(x(t)) (1+\lambda)[z(\varphi(t - \tau)) - z(x(t))]\;,\;
\lambda\in [-1,1]\}\;.
\]
Hence, if one finds a bound for \be\label{ofb} \ba{l} - W (x(t)) +
L_g V(x(t))
(1+\lambda)[z(\varphi(t - \tau)) - z(x(t))]\\
\hspace{5.5cm} + L_g V(x(t)) (\lambda - 1) z(\varphi(t -
\tau))\varrho\;, \ea\ee
then one also finds a bound for $v_a$.\\
Now,  inequality (\ref{ineq.h3'}) in Assumption (A3)  implies that
\be\label{khf} - \dst\frac{1}{2}W (x(t)) + L_g V(x(t))
(1+\lambda\delta)[z(\varphi(t - \tau)) - z(x(t))] \le
\kappa_{3}(R) [V(x(t)) + 1] \ee and that, for all $\lambda_a \in
[-1,1]$, \be\label{jkhf} - \dst\frac{1}{2}W (x) + L_g V(x)
(1+\lambda_a)[- z(x)] \le \kappa_{3}(R) [V(x) + 1] \ . \ee
Therefore, for all $\lambda \in [- 1,1]$, \be\label{jkh2f} -
\dst\frac{1}{2}W (x(t)) + L_g V(x(t)) (- 1 + \lambda) z(x(t))
\varrho\le \kappa_{3}(R) [V(x(t)) + 1] \ee with $\varrho = 1$ or
$0$. Next, from (\ref{ofb}), (\ref{khf}) and (\ref{jkhf}), we
deduce that
\[
\ba{rcl} && - W (x(t)) + L_g V(x(t))
(1+\lambda)[z(\varphi(t - \tau)) - z(x(t))]\\
&& \hspace{5.5cm}
+ L_g V(x(t)) (\lambda - 1) z(\varphi(t - \tau))\varrho\\
&\le & 2\kappa_{3}(R) [V(x(t)) + 1] =\kappa_{4}(R) [V(x(t)) + 1]
\;, \ea
\]
and therefore
$$
v_a \leq \kappa_{4}(R) [V(x(t)) + 1]\;.
$$
We deduce that necessarily, for all $t \in [0,T)$ such that $t\in
[t_i,t_{i+1})$, $i=0,1,\ldots,k$, (i.e.\ at all inter-switching
times) we have \be\label{yc} \ba{rcl} \dot V(t) & \leq &
\kappa_{4}(R) [V(x(t)) + 1] \ . \ea\ee On the other hand, for any
$t_i\in[0,T)$, with $i=0,1,\ldots,k$, (i.e.\ at the switching
times) \be\label{yc2} V(x(t_i^+)) = V(x(t_i)) \; . \ee We conclude
as in \cite{mazenc.bliman.tac06} that no finite escape time can
exist. Indeed, for any $t\in [0,T)$, let $t\in [t_i,t_{i+1})$ for
some $i$. Then, integrating (\ref{yc}) from $t_i$ to $t$, we
obtain
\[
V(x(t))+1\le {\rm e}^{\kappa_{4}(R)(t-t_i)}(V(x(t_i^+))+1)={\rm
e}^{\kappa_{4}(R)(t-t_i)}(V(x(t_i))+1)\;,
\]
where the latter equality follows from (\ref{yc2}). Similarly
\[\ba{rcl}
V(x(t_i))+1&\le &{\rm
e}^{\kappa_{4}(R)(t_i-t_{i-1})}(V(x(t_{i-1}^+))+1)\\
&= &{\rm
e}^{\kappa_{4}(R)(t_i-t_{i-1})}(V(x(t_{i-1}))+1)\\
&\le &{\rm e}^{\kappa_{4}(R)(t_i-t_{i-1})}
{\rm e}^{\kappa_{4}(R)(t_{i-1}-t_{i-2})}(V(x(t_{i-2}^+))+1)\\
&= &{\rm e}^{\kappa_{4}(R)(t_i-t_{i-1})}
{\rm e}^{\kappa_{4}(R)(t_{i-1}-t_{i-2})}(V(x(t_{i-2}))+1)\\
&\vdots&\\
&\le &{\rm e}^{\kappa_{4}(R)(t_i-t_{i-1})} {\rm
e}^{\kappa_{4}(R)(t_{i-1}-t_{i-2})}\ldots {\rm
e}^{\kappa_{4}(R)(t_{1}-t_{0})}
(V(x(t_{0}))+1)\\
&\le &{\rm e}^{\kappa_{4}(R)(t_i-t_{0})} (V(x(t_{0}))+1)\;. \ea\]
Recalling that $t_0=0$, it follows that:
\[
V(x(t))+1\le {\rm e}^{\kappa_{4}(R)t}(V(x(0))+1)\;,
\]
which shows that no finite escape time can actually exist. This
fact and \be\label{epi} |z(\varphi(t - \tau))| < (1 -
\delta)^{-1}u_0(R) \; , \; \forall t \in [0,\tau] \ee imply that
$x(t)$ can be extended beyond $T$. This yields a contradiction
with the definition of $T$. It follows that $x(t)$ is defined for
all $t\in [0, \tau]$. As before, by integrating (\ref{yc}) and
bearing in mind (\ref{yc2}), we infer that, for all $t \in
[0,\tau]$, \be\label{dtv3}\ba{rcl} V(x(t)) + 1 & \leq &
e^{\kappa_4(R) t} [V(x(0)) + 1] \; . \ea\ee It follows immediately
from (\ref{nom.stability}) that, for all $t \in [0,\tau]$,
\be\label{dtv4}\ba{rcl} \kappa_1(|x(t)|) & \leq & e^{\kappa_4(R)
\tau} [\kappa_2(|x(0)|) + 1] - 1
\\
& \leq & e^{\kappa_4(R) \tau} [\kappa_2(R) + 1] - 1 \; . \ea\ee It
follows that, for all $t \in [0,\tau]$, \be\label{dtv5}\ba{rcl}
|x(t)| & \leq & \alpha(R)\;, \ea\ee where $\alpha(\cdot)$ is the
function defined in (\ref{und}). Observe that (\ref{dtv5}) and the
inequality $||\varphi||_c \leq R$ imply that, for all $t \in
[\tau,2\tau]$, $|z(x(t - \tau))| \leq \sup_{|a|\le \alpha(R)}
|z(a)|$. Since $\omega(\cdot) \geq \alpha(\cdot)$, it follows
that, for all $t \in [\tau,2\tau]$, \be\label{dt7}\ba{rcl} |z(x(t
- \tau))| & < & (1-\delta)^{-1} u_0(R) \ . \ea\ee Hence, $\Psi
(z(x(t - \tau)))$ is well-defined for all $t\in [\tau, 2\tau]$.
Moreover, the time derivative of $z(x(t - \tau))$, namely
\[
\dst\frac{d}{dt}z(x(t - \tau)) = \left.\dst\frac{\partial
z(x)}{\partial x}\right|_{x=x(t - \tau)}[f(x(t - \tau))+g(x(t -
\tau))\Psi (z(\varphi(t-2\tau)))]\;,
\]
is bounded for all $t\in [\tau,2\tau]$, hence the length of the
inter-switching intervals is bounded away from zero on
$[\tau,2\tau]$, and therefore the switching times in that interval
do not accumulate in finite time. Next, arguing exactly as before
one can prove that $x(t)$ is defined for all $t\in [\tau, 2\tau]$
and that, for all $t \in [\tau,2\tau]$, \be\label{dtv8}\ba{rcl}
|x(t)| & \leq & \alpha(\alpha(R)) = \gamma(R)\;. \ea\ee
\subsection{Extending solutions for $t > 2\tau$}
To extend further the solution, we proceed by contradiction.
\\
Let us assume that \be\label{fca} \sup\left\{t : x(s) \; \;
\mbox{exists and} \; \; |x(s)| \leq \omega(R) \; , \; \forall \; s
\in [- 2 \tau, t]\right\} \ee is a finite real number that we
denote again $T$. From the inequality $||\varphi||_c \leq R$,
(\ref{dtv5}) and (\ref{dtv8}) and the facts that $R > 0$ and $W$
is positive definite, we deduce that $T > 2\tau$. Next, observe
that the continuity of the solutions and the definition of
$u_0(\cdot)$ in (\ref{u02}) and (\ref{ozbis}) imply that, for all
$t \in [- 2 \tau, T)$, \be\label{ras} |z(x(t))| \leq
\displaystyle\sup_{|a| \leq \omega(R)} \left\{|z(a)|\right\} <
\frac{u_0(R)}{1 - \delta} \ . \ee We exploit this inequality to
derive first an upper bound for $\dot V(t)$ and later on for $\dot
{\cal U}(t)$, where ${\cal U}$ is a Lyapunov-Krasowskii functional
to be introduced below.\\
Arguing as before (see (\ref{va}) and the sentence following it),
we claim that to find an upper bound for $\dot V(t)$, we need to
find an upper bound to the expression below for all $t \in [2
\tau, T)$, \be\label{oic} \ba{l} v_b(t)  := L_f V(x(t)) +
L_g V(x(t))(1 + \lambda\delta) z(x(t - \tau))\\
\hspace{4cm} + L_g V(x(t))(\lambda - 1) z(x(t - \tau)) \varrho
\ea\ee with $\lambda$ any number in the interval $[-1,1]$ and
$\varrho\in \{0,1\}$. Recall that  $\varrho=1$ if and only if
$|z(x(t - \tau))|\le (1+\delta)^{-1} u_j$. Thanks to Assumption
(A1), we deduce that, for all $t \in [2 \tau, T)$, \be\label{okc}
\ba{rcl} v_b(t) & \leq & - W (x(t)) + L_g V(x(t))(1 +
\lambda\delta) [z(x(t - \tau)) - z(x(t))]
\\
& & + L_g V(x(t))(\lambda - 1) z(x(t - \tau)) \varrho \;. \ea\ee
We now set $z(x(t-\tau))-z(x(t))$ in a form which allows us to use
Assumption (A2). Let $t_{i_j}$, $j=1,\ldots, k$ be the switching
times in the interval $(t-\tau, t)$ and set without loss of
generality $t_{i_0}=t-\tau$, $t_{i_{k+1}}=t$. We observe as before
that the switching times do not accumulate in finite time. Hence
we can write
\[\ba{l}
z(x(t)) - z(x(t-\tau))\\
= z(x(t)) - z(x(t_{i_k}))+ z(x(t_{i_k}))\ldots - z(x(t_{i_1}))+
z(x(t_{i_1}))
- z(x(t-\tau))\\
= \sum_{j=0}^{k} [z(x(t_{i_{j+1}})) - z(x(t_{i_{j}}))]\;. \ea\]
For each $m \in [t_{i_{j}}, t_{i_{j+1}})$, $j=0,1,\ldots,k$,
\[\ba{l}
z(x(t_{i_{j+1}}))
- z(x(t_{i_{j}}))\\
=\dst\int_{t_{i_{j}}}^{t_{i_{j+1}}} \dst\frac{\partial z}{\partial
x}(x(m))[f(x(m))+g(x(m))\Psi (z(x((t_{i_{j}}-\tau)^+)))]dm\\
= \dst\int_{t_{i_{j}}}^{t_{i_{j+1}}} \dst\frac{\partial
z}{\partial
x}(x(m))[f(x(m))+g(x(m))\Psi (z(x(m-\tau)))]dm\\
= \dst\int_{t_{i_{j}}}^{t_{i_{j+1}}} \dst\frac{\partial
z}{\partial x}(x(m))[f(x(m))+g(x(m))\psi_1(m-\tau)z(x(m-\tau))]dm
\ea\] where \footnote{ Let $\Psi (z(x(m-\tau)))=\tilde u_i$ for
$m\in [t_{i_{j}},t_{i_{j+1}})$, where $\tilde u_i=u_i$ or $\tilde
u_i=u_i(1+\delta)^{-1}$. Then $(1+\delta)^{-1}\tilde u_i <
z(x(m-\tau))\le (1-\delta)^{-1}\tilde u_i$. At each $m\in
[t_{i_{j}},t_{i_{j+1}})$,
$z(x(m-\tau))=\alpha(m-\tau)(1+\delta)^{-1}\tilde u_i+
(1-\alpha(m-\tau))(1-\delta)^{-1}\tilde u_i$, where
$\alpha(\cdot)$ takes value in $[0,1]$. Observe also that, since
$(1+\delta)^{-1}\tilde u_i$, $(1-\delta)^{-1}\tilde u_i$ are
constants, and $z(x(m-\tau))$ is a class $C^1$ function on the
interval $[t_{i_{j}},t_{i_{j+1}})$, so is $\alpha(m-\tau)$.
Moreover,
\[
\Psi (z(x(m-\tau)))=\tilde u_i=[\alpha(m-\tau)(1+\delta)^{-1}+
(1-\alpha(m-\tau))(1-\delta)^{-1}]^{-1}z(x(m-\tau))\;,
\]
where the function $\psi_1(m-\tau)=[\alpha(m-\tau)(1+\delta)^{-1}+
(1-\alpha(m-\tau))(1-\delta)^{-1}]^{-1}$ is a class $C^1$ function
which spans the interval $[1-\delta,1+\delta]$. Similar
considerations hold when $\tilde u_i$ is negative or equal to
zero. } $\psi_1(\cdot)$ is a class $C^1$ function taking value in
$[1-\delta,1+\delta]$. Overall we have
\[\ba{l}
z(x(t)) - z(x(t-\tau))\\
=\dst\int_{t-\tau}^{t} \dst\frac{\partial z}{\partial
x}(x(m))[f(x(m))+g(x(m))\psi_2(m-\tau)z(x(m-\tau))]dm \ea\] with
$\psi_2(\cdot)\,:\,[t-2\tau,t-\tau]\to [1-\delta,1+\delta]$ a
class $\overline C^0$ function, and hence
\[\ba{l}
v_b(t)  \leq  - W (x(t)) + L_g V(x(t)) (1 + \lambda\delta)
\dst\int_{t-\tau}^{t} \dst\frac{\partial z}{\partial
x}(x(m))[f(x(m))
\\
+ g(x(m)) \psi_2(m-\tau)z(x(m-\tau))]dm
 + \dst\frac{\partial V}{\partial x}g(x(t))
(\lambda - 1) z(x(t - \tau)) \varrho \;. \ea\] We deduce from
Assumption (A2), that we have, for all $t \in [2\tau,T)$,
\be\label{vjb} \ba{rcl} v_b(t) & \le & -\dst\frac{3}{4}W (x(t)) +
\dst\frac{1}{\Omega} \dst\int_{t-2\tau}^{t} W (x(\ell)) d\ell
\\
& & + L_g V(x(t)) (\lambda - 1) z(x(t - \tau)) \varrho \; . \ea\ee
Since $\varrho=1$ if and only if $|z(x(t - \tau))|\le
(1+\delta)^{-1} u_j$, we have that, for all $t \in [2\tau,T)$,
\[
\ba{rcl} v_b(t)  & \le & -\dst\frac{3}{4}W (x(t)) +
\dst\frac{1}{\Omega} \dst\int_{t-2\tau}^{t} W (x(\ell)) d\ell + 2
\left|L_g V(x(t))\right| (1+\delta)^{-1} u_j \;, \ea\] and
therefore \be\label{vjx} \ba{rcl} \dot V(t)  & \le &
-\dst\frac{3}{4}W (x(t)) + \dst\frac{1}{\Omega}
\dst\int_{t-2\tau}^{t} W (x(\ell)) d\ell + 2 \left|L_g
V(x(t))\right| (1+\delta)^{-1} u_j \;, \ea\ee
Next, with an abuse of notation, we define the following
Lyapunov-Krasowskii functional
\begin{equation}\label{U.bound}
\begin{array}{rcl}
{\cal U}(t) & = & V(x(t)) + \dst\frac{1}{8\tau}
\displaystyle\int_{t - 2 \tau}^{t} \displaystyle\int_{s}^{t}
W(x(\ell)) d\ell ds \ .
\end{array}
\end{equation}
We deduce from (\ref{vjx}) that, for all $t \in [2\tau,T)$, the
derivative of ${\cal U}$ along the trajectories of the system we
consider satisfies \be \ba{rcl} \dot {\cal U}(t) & \le & -
\dst\frac{1}{2}W (x(t)) + \left(\frac{1}{\Omega} -
\frac{1}{8\tau}\right)\dst\int_{t-2\tau}^{t} W (x(\ell)) d\ell +
2\left|L_g V(x(t))\right| (1+\delta)^{-1} u_j
\\
& \le & - \dst\frac{1}{2}W (x(t)) -
\frac{1}{16\tau}\dst\int_{t-2\tau}^{t} W (x(\ell)) d\ell + 2
\left|L_g V(x(t))\right| (1+\delta)^{-1} u_j \;, \ea\ee where the
last inequality is a consequence of the condition $\Omega \geq 16
\tau$. Let $\kappa_5$ be a positive increasing function of class
$C^1$ such that for all $x \in \R^n$ \be\label{N} \left|L_g
V(x)\right| \le \kappa_5(V(x)), \ee and set $\kappa_6(R) =
\kappa_5(\kappa_2(\omega(R)))$ (observe that $\kappa_6$ is
continuous). From (\ref{N}), (\ref{nom.stability}), the definition
of $\kappa_6$ and the definition of $T$, we infer that, for all $t
\in [2\tau,T)$, \be\label{ejf} \ba{rcl} \dot {\cal U}(t) & \le & -
\dst\frac{1}{2}W (x(t)) - \frac{1}{16\tau}\dst\int_{t-2\tau}^{t} W
(x(\ell)) d\ell + 2\kappa_6(R) (1+\delta)^{-1} u_j
\\
& \le & - \dst\frac{1}{2}W (x(t)) -
\frac{1}{32\tau^2}\dst\int_{t-2\tau}^{t}\dst\int_{s}^{t} W
(x(\ell)) d\ell ds + 2\kappa_6(R) (1+\delta)^{-1} u_j \; . \ea\ee
Next, we would like to express the first two terms on the
right-hand side of the last inequality in terms of ${\cal U}(t)$.
This is possible according to Lemma \ref{inap} in Appendix
\ref{appen}. Namely, one can determine a $C^1$ class-${\cal
K}_{\infty}$ function $\kappa_{7\tau}$ and a function
$\kappa_{8\tau}$ of class $C^1$, positive and nondecreasing such
that, for all $x\in \R^n$ and $z \geq 0$, \be\label{M}
\frac{\kappa_{7\tau}\left(V(x) +
z\right)}{\kappa_{8\tau}\left(V(x) + z\right)} \le
\dst\frac{1}{2}W(x) + \dst\frac{1}{4\tau} z \;, \ee where
$z=\frac{1}{8\tau}\dst\int_{t-2\tau}^{t}\dst\int_{s}^{t} W
(x(\ell)) d\ell ds$. From (\ref{M}), it is possible to deduce
that, for all $t \in [2\tau,T)$, \be\label{gdf} \ba{rcl} \dot
{\cal U}(t) & \le & - \dst\frac{\kappa_{7\tau}\left({\cal
U}(t)\right)}{\kappa_{8\tau}\left({\cal U}(t)\right)} +
2\kappa_6(R) (1+\delta)^{-1} u_j \; . \ea\ee Since, for all $t \in
[2\tau,T)$, $|x(t)| \leq \omega(R)$, we deduce that for all $t \in
[2\tau,T)$, $V(x(t)) \leq \kappa_2(\omega(R))$ and
$\dst\int_{t-2\tau}^{t}\dst\int_{s}^{t} W (x(\ell)) d\ell ds \leq
2\tau^2 \displaystyle\sup_{|a| \leq \omega(R)} W (a)$. It follows
that \be\label{gdc} \ba{rcl} {\cal U}(t) & \le &
\kappa_2(\omega(R))+ 2\tau^2 \displaystyle\sup_{|a| \leq
\omega(R)} W (a) =: \kappa_{9}(R) \ea\ee where $\kappa_{9}$ is
continuous and nondecreasing. Next, let us prove that, for all $t
\in [2\tau,T)$, the inequality \be \label{beb} {\cal U}(t) \leq
\kappa_2(\gamma(R)) + \frac{\tau}{4}\displaystyle\sup_{|a| \leq
\gamma(R)} W (a) \ee is satisfied. This result is the consequence
of (\ref{gdf}) and the fact that ${\cal U}(2\tau) \leq
\kappa_2(\gamma(R)) + \frac{\tau}{4}\displaystyle\sup_{|a| \leq
\gamma(R)} W (a)$ and \be\label{pli} \ba{rcl} \dot {\cal U}(t) & <
& 0 \ea\ee when ${\cal U}(t) = \kappa_2(\gamma(R)) +
\frac{\tau}{4}\displaystyle\sup_{|a| \leq \gamma(R)} W (a)$
{\em provided that} $u_j$ is appropriately chosen.\\
To see this, observe  in particular (recall (\ref{yeio})) that
$u_j$ can be made small by increasing the number of quantization
levels $j$. Namely, let \be\label{j02} j \geq
\left|\left(\log\dst\frac{\mu(\varepsilon,R)(1 +
\delta)}{u_0(R)}\right) \left(\log\dst\frac{1 - \delta}{1 +
\delta}\right)^{-1} \right| + 1 \, \ee where $\mu$ is continuous
and such that, for all the real numbers $\varepsilon > 0, R > 0$,
\be\label{ryt} 0 < \mu(\varepsilon,R) \leq \min\left\{A_1(R),
A_2(\varepsilon,R)\right\} \ee with \be\label{yak} \ba{rcl} A_1(R)
& = & \dst\frac{\kappa_{7\tau}\left(\kappa_2(\gamma(R)) + \dst
\frac{\tau}{4} \displaystyle\sup_{|a| \leq \gamma(R)}
W(a)\right)}{4 \kappa_6(R) \kappa_{8\tau}\left(\kappa_2(\gamma(R))
+ \dst \frac{\tau}{4} \displaystyle\sup_{|a| \leq \gamma(R)}
W(a)\right)} \; ,
\\
A_2(\varepsilon,R) & = &
\dst\frac{\kappa_{7\tau}(\kappa_1(\varepsilon))}{4\kappa_6(R)
\kappa_{8\tau}\left(\kappa_2(\gamma(R)) + \dst \frac{\tau}{4}
\displaystyle\sup_{|a| \leq \gamma(R)} W(a)\right)} \; . \ea\ee
Observe that (\ref{j02}) and  (\ref{yeio})  imply \be\label{mu}
(1+\delta)^{-1}u_{j} \le \mu(\varepsilon,R) \ . \ee Hence,
(\ref{gdf}) rewrites as \be\label{gdf.bis} \ba{rcl} \dot {\cal
U}(t) & \le & - \dst\frac{\kappa_{7\tau}\left({\cal
U}(t)\right)}{\kappa_{8\tau}\left({\cal U}(t)\right)} +
2\kappa_6(R) \mu(\varepsilon,R) \; . \ea\ee From this, since
$\mu(\varepsilon,R)\le A_1(\varepsilon,R)$, it is immediate to see
that when ${\cal U}(t) = \kappa_2(\gamma(R)) +
\frac{\tau}{4}\displaystyle\sup_{|a| \leq \gamma(R)} W (a)$,
(\ref{pli}) is satisfied and (\ref{beb}) holds.
\\
From (\ref{beb}), it follows immediately that, for all $t \in
[2\tau,T)$, \be\label{pg2i} |x(t)| \leq
\kappa_1^{-1}\left(\kappa_2(\gamma(R)) +
\frac{\tau}{4}\displaystyle\sup_{|a| \leq \gamma(R)} W (a)\right)
\ . \ee We deduce from (\ref{berl}) and $R > 0$ that, for all $t
\in [2\tau,T)$, \be\label{f1} |x(t)| < \omega(R) \ . \ee This
inequality and (\ref{ras}) imply that $x(t)$ can be extended
beyond $T$. This yields a contradiction with the definition of
$T$. We deduce that $x(t)$ is defined over $[- 2 \tau, + \infty)$
and bounded in norm by $\omega(R)$.
\subsection{Practical convergence}
Observe that, arguing as before, one can prove that for all $t
\geq 2\tau$, we have \be\label{gdb1} \ba{rcl} \dot {\cal U}(t) &
\le & - \dst\frac{\kappa_{7\tau}\left({\cal
U}(t)\right)}{\kappa_{8\tau}\left({\cal U}(t)\right)} +
2\mu(\varepsilon,R) \kappa_6(R)
\\
& \le & - \frac{\kappa_{7\tau}\left({\cal
U}(t)\right)}{\kappa_{8\tau}\left( \kappa_2(\gamma(R)) +
\frac{\tau}{4} \displaystyle\sup_{|a| \leq \gamma(R)} W (a)
\right)} + 2\mu(\varepsilon,R) \kappa_6(R) \; . \ea\ee Since,
according to (\ref{ryt}), \be\label{gdb2} \ba{rcl}
4\mu(\varepsilon,R) \kappa_6(R) & \le & \dst
\frac{\kappa_{7\tau}(\kappa_1(\varepsilon))}{\kappa_{8\tau}\left(
\kappa_2(\gamma(R)) + \frac{\tau}{4} \displaystyle\sup_{|a| \leq
\gamma(R)} W (a) \right)} \ea\ee we deduce that, there exists $t_L
\geq 0$ such that, for all $t \geq t_L$, the inequality
\be\label{gdb3} \ba{rcl} {\cal U}(t) & \le & \kappa_1(\varepsilon)
\ea\ee is satisfied. It follows that, for all $t \geq t_L$,
\be\label{gd2} \ba{rcl} |x(t)| & \le & \varepsilon \ , \ea\ee that
is the thesis.

\begin{remark}
It has been observed in the proof that, by the definition
(\ref{fca}) of $T$, for all $t\in [2\tau,T)$, $|x(t)|<\omega(R)$,
and therefore
$\frac{1}{8\tau}\dst\int_{t-2\tau}^{t}\dst\int_{s}^{t} W (x(\ell))
d\ell ds \leq \frac{\tau}{4} \sup_{|a| \leq \omega(R)} W (a)$.
Set:
\begin{equation}
\label{ells} \kappa_{7\tau}^R(\xi) =
\frac{\kappa_{7\tau}(\xi)}{\kappa_{8\tau}\left(\kappa_2(\omega(R))
+ \displaystyle\frac{\tau}{4}\sup_{|a| \leq \omega(R)} W
(a)\right)} \ .
\end{equation}
Then, one can use in the proof the inequality \be\label{gdf.bis}
\ba{rcl} \dot {\cal U}(t) & \le & - \kappa^R_{7\tau}\left({\cal
U}(t)\right) + 2\mu(\varepsilon,R) \kappa_6(R) \; , \ea\ee instead
of (\ref{gdf}). In particular one can follow exactly the same
passages as before, provided that in the definition (\ref{ryt}) of
$\mu(\varepsilon,R)$, the functions $A_1(R)$ and
$A_2(\varepsilon,R)$ in (\ref{yak}) are defined as
\be\label{yak.bis} \ba{rcl} A_1(R) & = &
\dst\frac{\kappa^R_{7\tau}\left(\kappa_2(\gamma(R)) + \dst
\frac{\tau}{4} \displaystyle\sup_{|a| \leq \gamma(R)}
W(a)\right)}{4 \kappa_6(R)} \; ,
\\
A_2(\varepsilon,R) & = &
\dst\frac{\kappa^R_{7\tau}(\kappa_1(\varepsilon))}{4\kappa_6(R)}
\;. \ea\ee Then, by replacing the differential inequality
(\ref{gdb1}) with (\ref{gdf.bis}), and the inequality (\ref{gdb2})
with \be\label{gdb2.bis} \ba{rcl} 4\mu(\varepsilon,R) \kappa_6(R)
& \le & \kappa_{7\tau}^R(\kappa_1(\varepsilon))\;, \ea\ee we can
again conclude that $x(t)$ enters the closed ball of radius
$\varepsilon$
in finite time and remains in it thereafter.\\
This remark is useful to simplify the proof in the particular case
where a constant function can be chosen for the function
$\kappa_{8\tau}$ in (\ref{M}), for instance when the positive
definite function $W (x)$ is lower bounded by a class-${\cal
K}_\infty$
function, as it happens when system (\ref{system.f.ol}) is linear.\\
Observe, finally, that a function $\kappa_{7\tau}^R$ of class
${\cal K}$ such that (\ref{gdf.bis}) is satisfied can be found
without necessarily relying on the knowledge of $\kappa_{7\tau}$
and $\kappa_{8\tau}$. In fact, bearing in mind (\ref{ejf}), it
suffices to find $\kappa_{7\tau}^R$ such that
\[
\dst\frac{1}{2}W (x(t)) + \frac{1}{4\tau} z(t)\ge
\kappa_{7\tau}^R({\cal U}(t))
\]
with
\[
z(t)=\frac{1}{8\tau}\dst\int_{t-2\tau}^{t}\dst\int_{s}^{t} W
(x(\ell)) d\ell ds\;.
\]
For instance, one can choose
$\kappa_{7\tau}^R(s)=K_e(\frac{1}{2}B_S(s))$, where
\[
K_e(m)=\dst\frac{1}{\overline m}\dst\int_0^m \min_{l\le |\xi|\le
\overline m} {\cal W}(\xi)dl\;,\quad m\in [0, \overline m]\;,
\]
${\cal W}(\xi)=\frac{1}{2} W (x) + \frac{1}{4 \tau} z$,
$\xi=(x^T\, z)^T$, $\overline m=\omega(R)+z_R$ and $B_S$ is
defined as \be \label{rt2e} B_S(l) =
\min\left\{\kappa_2^{-1}\left(\frac{l}{2}\right),
\frac{l}{2}\right\} \;. \ee As a matter of fact, for all $|\xi|\le
\overline m$,
\[\ba{rcl}
{\cal W}(\xi) &\ge & \dst\min_{|\xi|\le \overline m} {\cal W}(\xi)
\ge \dst\frac{1}{\bar m}\dst\int_0^{|\xi|} \min_{|\xi|\le
\overline m} {\cal W}(\xi) dl \ge \dst\frac{1}{\bar
m}\dst\int_0^{|\xi|} \min_{l\le |\xi|\le \overline m} {\cal
W}(\xi) dl =K_e(|\xi|)\;. \ea\] Now, $|x|+z\le 2|\xi|$ and \be
\label{sre} |x| + z \geq \kappa_2^{-1}(V(x)) + z \geq B_S(V(x) +
z) \ee with \footnote{ Let $\alpha_1(r)=\kappa_2^{-1}(r)$,
$\alpha_2(r)=r$, $\alpha_3(r)=\min\{\alpha_1(r),\alpha_2(r)\}$,
and $2a=V(x)$, $2b=z$. Then, bearing in mind that, for any
function $\alpha$ of class ${\cal K}_\infty$, $\alpha(a+b)\le
\alpha(2a)+\alpha(2b)$, we have
\[
\alpha_1(2a)+\alpha_2(2b)\ge \alpha_3(2a)+\alpha_3(2b)\ge
\alpha_3(a+b)\;.
\]
This proves $\kappa_2^{-1}(V(x)) + z \geq B_S(V(x) + z)$. }
$B_S(l)$ as in (\ref{rt2e}). Hence,
\[
{\cal W}(\xi)\ge K_e(|\xi|)\ge K_e(\frac{1}{2}(|x| + z)) \ge
K_e(\frac{1}{2}B_S(V(x) + z))=\kappa_{7\tau}^R(V(x) + z)
\]
and therefore
\[
\frac{1}{2} W (x(t)) + \frac{1}{4 \tau} z(t)={\cal W}(\xi(t))\ge
\kappa_{7\tau}^R(V(x(t)) + z(t))= \kappa_{7\tau}^R({U}(t))
\]
as desired.\\
We could have stated the   result directly in terms of the
class-${\cal K}$ function $\kappa_{7\tau}^R$ just derived rather
than introducing the two class-${\cal K}_\infty$ functions
$\kappa_{7\tau},\kappa_{8\tau}$. We decided to adopt the latter in
order not to have in $A_1(R), A_2(\varepsilon,R)$ (and hence in
the conditions on the number of quantization levels $j$) a
class-${\cal K}$ function depending {\em implicitly} on the
parameter $R$.
\end{remark}

\begin{remark}
The proof is constructive: It gives the explicit expressions for
the two design parameters $u_0$ (see (\ref{u02})), and $j$ (see
(\ref{j02}) and (\ref{ryt})).
\end{remark}

\section{Conclusion}
We have presented a Lyapunov-Krasowskii functional approach to
solve the problem of determining quantized feedbacks with delay
which semi-globally practically stabilize the origin of nonlinear
systems. For a fairly general family of systems, and given any
value of the quantization density, we have characterized the
maximal allowable constant delay which the closed-loop system can
tolerate. A problem which in our opinion would be interesting to
investigate is how, for systems with a well-defined relative
degree, our result can be propagated via the backstepping
technique.

\appendix

\section{Technical lemmas}
\label{appen}

\begin{lemma}
\label{iznap} Let ${\cal W} : \mathbb{R}^n \rightarrow \mathbb{R}$
be a continuous and positive definite function. For all $m \geq
0$, let \be\label{ipa} \ba{rcl} K_a(m) & = &
\displaystyle\int_{0}^{\min\{m,1\}} \displaystyle\min_{l \leq
|\xi| \leq 1} {\cal W}(\xi) dl + \max\{0,m - 1\} ,
\\
K_b(m) & = & 1 + \frac{K_a(m)}{\displaystyle\min_{1 \leq |\xi|
\leq \max\{1,m\}} {\cal W}(\xi)} \ . \ea\ee Then $K_a$ belongs to
${\cal K}_\infty$, $K_b$ is continuous, positive and increasing
over $[0,+\infty)$ and, for all $X \in \mathbb{R}^n$, \be
\label{iza} \frac{K_a(|X|)}{K_b(|X|)} \leq {\cal W}(X) \ . \ee
\end{lemma}

\noindent {\bf Proof.} The fact that ${\cal W}(X)$ is positive
definite implies that both $K_a$ and $K_b$ are well-defined and
continuous. Let us prove that $K_a$ belongs to ${\cal K}_\infty$.
Observe that $K_a(0) = 0$. When $m \in [0,1]$, $K_a(m) =
\displaystyle\int_{0}^{m} \displaystyle\min_{l \leq |\xi| \leq 1}
{\cal W}(\xi) dl$. Therefore this function is increasing over
$[0,1]$. When $m > 1$, $K_a(m) = \displaystyle\int_{0}^{1}
\displaystyle\min_{l \leq |\xi| \leq 1} {\cal W}(\xi) dl + m - 1$.
Therefore this function is increasing over $[1,+\infty)$ and goes
to the infinity when its argument does. Consequently, $K_a$ is of
class ${\cal K}_\infty$. If follows that $K_b$ is a positive and
increasing over
$[0,+\infty)$.\\
Next, to establish (\ref{iza}), we distinguish between two cases. \\
{\em First case:} $|X| \leq 1$. Then $K_a(|X|) =
\displaystyle\int_{0}^{|X|} \displaystyle\min_{l \leq |\xi| \leq
1} {\cal W}(\xi) dl \leq |X| \displaystyle\min_{|X|\leq |\xi| \leq
1} {\cal W}(\xi) \leq {\cal W}(X)$.
Moreover, $K_b(|X|) \geq 1$. It follows that $\frac{K_a(|X|)}{K_b(|X|)} \leq {\cal W}(X)$. \\
{\em Second case:} $|X| \geq 1$. Then \be\label{ypa} \ba{rcl}
K_a(|X|) & = & \displaystyle\int_{0}^{1} \displaystyle\min_{l \leq
|\xi| \leq 1} {\cal W}(\xi) dl + |X| - 1 > 0\;,
\\
K_b(|X|) & = & 1 + \frac{K_a(|X|)}{\displaystyle\min_{1 \leq |\xi|
\leq |X|} {\cal W}(\xi)}
> \frac{K_a(|X|)}{\displaystyle\min_{1 \leq |\xi| \leq |X|} {\cal W}(\xi)} >
0\;. \ea\ee Therefore \be \ba{rcl} \displaystyle\min_{1 \leq |\xi|
\leq |X|} {\cal W}(\xi) & > & \frac{K_a(|X|)}{K_b(|X|)} > 0\;.
\ea\ee It follows that \be\label{ypb} {\cal W}(X) >
\frac{K_a(|X|)}{K_b(|X|)} > 0 \ . \ee

\begin{lemma}
\label{inap} Let $\tau > 0$, $W $ be a positive definite function.
Then one can determine a function $K_c$ of class ${\cal K}_\infty$
and a function $K_d$, positive, continuous and increasing over
$[0,+\infty)$ such that, for all $x \in \mathbb{R}^n$ and $z \geq
0$, \be \label{iya} \frac{K_c(V(x) + z)}{K_d(V(x) + z)} \leq
\frac{1}{2} W (x) + \frac{1}{4 \tau} z\;. \ee
\end{lemma}

\noindent {\bf Proof.} First, observe that the inequalities
(\ref{nom.stability}) imply that for all $x \in \mathbb{R}^n$, $z
\geq 0$, \be \label{tre} |x| + z \leq \kappa_1^{-1}(V(x)) + z \leq
B_L(V(x) + z) \ee with \be \label{rt1e} B_L(l) = \kappa_1^{-1}(l)
+ l \ee and
\[
|x| + z \geq \kappa_2^{-1}(V(x)) + z \geq B_S(V(x) + z)
\]
with
\[
B_S(l) = \min\left\{\kappa_2^{-1}\left(\frac{l}{2}\right),
\frac{l}{2}\right\} \; .
\]
From Lemma \ref{iznap}, it follows immediately that one can
determine a function $K_a$ of class ${\cal K}_\infty$ and a
function $K_b$, positive, continuous and increasing over
$[0,+\infty)$ such that, for all $x \in \mathbb{R}^n$ and $z \geq
0$, \be \frac{K_a(|x| + z)}{K_b(|x| + z)} \leq \frac{1}{2} W (x) +
\frac{1}{4 \tau} z \ . \ee From this inequality, (\ref{tre}) and
(\ref{sre}), it follows that \be \frac{K_a(B_S(V(x) +
z))}{K_b(B_L(V(x) + z))} \leq \frac{1}{2} W (x) + \frac{1}{4 \tau}
z \ . \ee This allows us to conclude.

\end{document}